\documentclass[numbib]{imamat}

\usepackage{amsmath}
\usepackage{amsthm}
\usepackage{cite}
\usepackage{hyperref}
\usepackage{graphics}
\usepackage{algorithmic}
\usepackage{subfig}
\usepackage{url}

\usepackage{tikz, subfig, amsthm, multirow, float, enumerate}
\usepackage{tikz-3dplot}
\usetikzlibrary{arrows, shapes, positioning}
\usepackage{bm}

\usepackage{amsthm}
\newtheoremstyle{theorem}{6pt}{6pt}{\rm}{}{\sffamily}{ }{ }{}
\theoremstyle{theorem}

\newtheoremstyle{lemma}{6pt}{6pt}{\rm}{}{\sffamily}{ }{ }{}
\theoremstyle{lemma}

\newtheoremstyle{example}{6pt}{6pt}{\rm}{}{\sffamily}{ }{ }{}
\theoremstyle{example}

\newtheoremstyle{corollary}{6pt}{6pt}{\rm}{}{\sffamily}{ }{ }{}
\theoremstyle{corollary}

\newtheoremstyle{definition}{6pt}{6pt}{\rm}{}{\sffamily}{ }{ }{}
\theoremstyle{definition}

\newtheoremstyle{remark}{6pt}{6pt}{\rm}{}{\sffamily}{ }{ }{}
\theoremstyle{remark}

\newtheoremstyle{approximation}{6pt}{6pt}{\rm}{}{\sffamily}{ }{ }{}
\theoremstyle{approximation}

\newtheoremstyle{scheme}{6pt}{6pt}{\rm}{}{\sffamily}{ }{ }{}
\theoremstyle{scheme}

\DeclareMathOperator*{\argmin}{arg\,min}
\DeclareMathOperator*{\argmax}{arg\,max}

\newcommand{\REP}{\mathrm{LREP}}
\newcommand{\DN}{\Delta_N}
\newcommand{\thetaidx}{q(N)}
\newcommand{\thetaN}{\boldsymbol \theta_{\thetaidx}}

\newcommand{\elt}{A_{N}(\thetaN) }
\newcommand{\Gam}{B_{N}(\thetaN) }

\newcommand{\Gamc}{C_{N}(\thetaN) }

\usepackage{color}
\newcommand{\ak}[1]{{\color{black} #1}}

\AtBeginDocument{
  \label{FirstPage2}
  \def\fpage{\pageref{FirstPage2}}
}

\let\BeginKnitrBlock\begin \let\EndKnitrBlock\end
\begin{document}

\title{On the instability and degeneracy of deep learning models}
\shorttitle{Instability of deep learning models}
\shortauthorlist{Kaplan, Nordman, and Vardeman}

\author{{\sc Andee Kaplan}\(^*\), \\[2pt] Department of Statistical Science, Duke University, P.O. Box 90251, Durham, NC 27708-0251, USA \\ \(^*\){\email{Corresponding author: \href{mailto:andrea.kaplan@duke.edu}{\nolinkurl{andrea.kaplan@duke.edu}}}}
\\[6pt] {\sc Daniel Nordman} \\[2pt] Department of Statistics, Iowa State University, 2438 Osborn Dr., Ames, IA 50011-1090, USA \\ {\href{mailto:dnordman@iastate.edu}{\nolinkurl{dnordman@iastate.edu}}}  
\\[6pt] {\sc Stephen Vardeman} \\[2pt] Departments of Statistics and Industrial and Manufacturing Systems Engineering, Iowa State University, 2438 Osborn Dr., Ames, IA 50011-1090, USA \\ {\href{mailto:vardeman@iastate.edu}{\nolinkurl{vardeman@iastate.edu}}}  
}

\maketitle

\begin{abstract}
{A probability model exhibits instability if small changes in a data outcome result in large, and often unanticipated, changes in probability. This instability is a property of the probability \ak{model, given by a distributional form and a given configuration of parameters}. For correlated data structures found in several application areas, there is increasing interest in \ak{identifying} such sensitivity in model probability structure. We consider the problem of quantifying instability for general probability models defined on sequences of observations, where each sequence of length \(N\) has a finite number of possible values \ak{that can be taken at each point}. A sequence of probability models results, indexed by \(N\), \ak{and an associated parameter sequence,} that accommodates data of expanding dimension. Model instability is formally shown to occur when a certain log-probability ratio under such models grows faster than \(N\). In this case, a one component change in the data sequence can shift probability by orders of magnitude. Also, as instability becomes more extreme, the resulting probability models are shown to tend to degeneracy, placing all their probability on potentially small portions of the sample space. These results on instability apply to large classes of models commonly used in random graphs, network analysis, and machine learning contexts.}
{Degeneracy, Instability, Classification, Deep Learning, Graphical Models}
\end{abstract}

\hypertarget{introduction}{%
\section{Introduction}\label{introduction}}

We consider the behavior, and the potential impropriety, of sequences of discrete probability models built to incorporate observations of increasing sample size \(N\). Interest is in identifying instability in such models, which is roughly characterized by probabilities with extreme sensitivity to small changes in data configuration. The concept of instability was introduced in the field of statistical physics (i.e., point processes) by \citet{ruelle1999statistical} and then further extended by \citet{schweinberger2011instability} for a family of exponential models. At issue, models exhibiting instability are typically undesirable as these tend to provide poor representations of data or data-generation. As an example, such models can include near-degenerate distributions that assign essentially all probability mass to only a subset of an overall sample space. The latter issue in connection to degeneracy has been recognized as a concern in that dominant model outcomes may not resemble observed data \citep[cf.~][]{handcock2003assessing}. As a compounding issue, model instability often has direct negative impacts for statistical inference and computations based on likelihood functions. Namely, volatilities in probability structure can potentially hamper the numerical evaluations required for maximum likelihood estimation as well as other model-based simulations via Markov Chain Monte Carlo (MCMC). These reasons motivate our general study of instability for a broad class of probability models, described next.

In the model framework, let \(\boldsymbol X_N = (X_1, \dots, X_N)\) denote a collection of discrete random variables with a finite sample space, \(\mathcal{X}^N\), represented as some \(N\)-fold Cartesian product. That is, \(\mathcal{X}\) with \(|\mathcal{X}| < \infty\) denotes the set of potential outcomes for each single variable \(X_i\), so that the product space \(\mathcal{X}^N\) corresponds to values for the variables \(\boldsymbol X_N=(X_1,\ldots,X_N)\). For each \(N\), let \(P_{\ak{\thetaN}}\) denote a probability model on \(\mathcal{X}^N\), under which \(P_{\ak{\thetaN}}(x_1, \dots, x_N) > 0\) is the probability of the data outcome \((x_1, \dots, x_N) \in \mathcal{X}^N\). In this, we assume that the model support of \(P_{\ak{\thetaN}}\) is the sample space \(\mathcal{X}^N\). This framework produces probability models \(P_{\ak{\thetaN}}\), indexed by a \ak{defining} sequence of parameters \(\ak{\thetaN}\), to describe data \(\boldsymbol X_N\) of any given sample size \(N \geq 1\). For simplicity, we will refer to this distributional class as \emph{Finite Outcome Everywhere Supported (FOES)} models in the following. The dimension and structure of such parameters are generic, without restriction, though natural cases will be seen to include those where \(\ak{\thetaN} \in \mathbb{R}^{q(N)}\) for some arbitrary integer-valued function \(q(\cdot) \geq 1\).

Section \ref{examples} provides some examples of FOES models encountered in graph/network analysis and machine learning (i.e., deep learning models). These are used as references for later illustrations. Section \ref{instability-results} then establishes several formal results for FOES models with regard to instability. \citet{schweinberger2011instability} originally developed instability results specific to a certain class of discrete exponential models. For similar exponential models with random networks, \citet{handcock2003assessing} studied model degeneracy, where a probability model places near complete mass on modes and may thereby narrow the feasible model outcomes. As findings here and from \citet{schweinberger2011instability} suggest, model instability and degeneracy may also be related by viewing degeneracy as an extreme, or limiting form, of instability. Our main results establish a broad characterization of model instability, appropriate across the whole FOES model class, that incorporates results of \citet{schweinberger2011instability} as a special case. We prescribe a general and simple condition for identifying instability in a FOES model sequence, which quantifies whether certain maximal probabilities in a FOES model are too extreme relative to the sample size \(N\). When these conditions are met, the probability structure of a FOES model is shown to exhibit extreme sensitivity, with probability assignments possessing extreme peaks and troughs across nearly identical outcomes. As the measure of model instability increases, probabilities from an unstable FOES model additionally increase in volatility and provably slide into degeneracy. Section \ref{implications} then emphasizes the implications of such model instability, showing that such impropriety can be expected to numerically hinder maximum likelihood estimation and MCMC-based simulations. As one potential remedy, suggestions are given for constraining model parameterizations to avoid the most problematic regions of the parameter space \ak{(see Section \ref{bayes} for further discussion)}. Proofs of the main results appear in Appendix \ref{appendix-instab}.

\hypertarget{examples}{%
\section{Examples}\label{examples}}

Many model families \ak{lead to} FOES models. \ak{(For clarity, we use ``model family" to refer to generic class or functional form of parametric distributions, while "model" often refers here to a member of the class for an instantiation or selection of parameters. This distinction is not so important for the examples presented in Sections \ref{discrete-exponential-family-models}-\ref{deep-learning}, but becomes more so in Section \ref{instability-results} where model instability depends on both form and a given parameter configuration)} For illustration, \ak{Sections \ref{discrete-exponential-family-models}-\ref{deep-learning} present} three specific examples of FOES model \ak{classes}, including models with deep architectures.

\hypertarget{discrete-exponential-family-models}{%
\subsection{Discrete Exponential Family Models}\label{discrete-exponential-family-models}}

For random variables \(\boldsymbol X \equiv\boldsymbol X_N= (X_1, \dots, X_N)\) with sample space \(\mathcal{X}^N\), \(|\mathcal{X}| < \infty\), consider an exponential \ak{model family} for \(\boldsymbol X\) with probability mass function given by
\begin{equation}
\label{eq:expo}
p_{N, \boldsymbol \theta}(\boldsymbol x) = \exp\left[\boldsymbol\eta^T(\boldsymbol \theta) \boldsymbol g_N(\boldsymbol x) - \psi(\boldsymbol \theta)\right], \qquad \boldsymbol x \in \mathcal{X}^N,
\end{equation}
depending on parameter vector \(\boldsymbol \theta \in \Theta_{\ak{\thetaidx}} \subset \mathbb{R}^{k}\) and natural parameter function \(\boldsymbol \eta : \mathbb{R}^k \mapsto \mathbb{R}^L\) with fixed positive integers \(k\) and \(L\) denoting their dimensions. Above, \(\boldsymbol g_N : \mathcal{X}^N \mapsto \mathbb{R}^L\) is a vector of sufficient statistics, while
\[
\psi(\boldsymbol \theta) = \log \sum\limits_{\boldsymbol x \in \mathcal{X}^N}\exp\left[\boldsymbol \eta^T(\boldsymbol \theta) \boldsymbol g_N(\boldsymbol x) \right], \qquad \boldsymbol \theta \in \Theta_\ak{\thetaidx}\equiv \{\boldsymbol \theta \in \mathbb{R}^k : \psi(\boldsymbol \theta) < \infty \},
\]
denotes the normalizing function with parameter space \(\Theta_{\ak{\thetaidx}}\). The natural parameter function \(\eta (\boldsymbol \theta)\) has a linear form (i.e., \(\eta (\boldsymbol \theta)= \bm{A} \boldsymbol \theta\) for a given \(L \times k\) matrix \(\bm{A}\)) in many common model formulations, though may also be nonlinear (e.g., curved exponential families). In the linear case, \(\eta (\boldsymbol \theta) = \boldsymbol \theta\) may be generally assumed in the exponential parameterization with a minor modification to the definition of sufficient statistics \(\boldsymbol g_N(\boldsymbol x)\).

Such discrete \ak{exponential models} are special cases of the FOES models, as seen by defining \(P_{\ak{\thetaN}}(\boldsymbol x)\equiv p_{N, \ak{\thetaN}}(\boldsymbol x)> 0\), \(\boldsymbol x \in \mathcal{X}^N\), based on \eqref{eq:expo} and a parameter sequence \(\ak{\thetaN} \in \Theta_{\ak{\thetaidx}} \subset \mathbb{R}^k\). For example, if observations \(\boldsymbol X = (X_1,\ldots,X_N)\) correspond to \(N\) independent and identically distributed Bernoulli random variables, each indicating a binary \(0\)-\(1\) outcome, the resulting probabilities have exponential form \eqref{eq:expo} given by
\begin{equation}
\label{eq:mod1}
P_{\ak{\thetaN}}(\boldsymbol x) \propto
 \exp\left[\ak{\thetaN} \sum_{i=1}^N x_i\right], \qquad \boldsymbol x=(x_1,\ldots,x_N) \in\{0,1\}^N, 
 \end{equation}
with sufficient statistic \(\boldsymbol g_N(\boldsymbol x)\equiv \sum_{i=1}^N x_i\) and ``log odds ratio'' parameter \(\ak{\thetaN} \equiv \log[ P_{\ak{\thetaN}}(X_i=1)/P_{\ak{\thetaN}}(X_i=0) ] \in \mathbb{R}\). More generally, supposing \(\boldsymbol X =(X_1,\ldots,X_N)\) represent \(N\) independent trials, each assuming an outcome \(\{1,\ldots,k\}\) among \(k\) possibilities (e.g., a die roll), a multinomial distribution is given by
\begin{equation}
\label{eq:mod11}
P_{\ak{\thetaN}}(\boldsymbol x) \propto  \exp\left[  \ak{\thetaN}^T g_N(\boldsymbol x)   \right] =
\exp\left[ \sum_{j=1}^k {\theta_{j,\ak{\thetaidx}}} \sum_{i=1}^N \mathbb{I}(x_i=j) \right], \qquad \boldsymbol x  \in\{1,\ldots,k\}^N, 
\end{equation}
with sufficient statistic \(\boldsymbol g_N(\boldsymbol x)\) involving a count \(\sum_{i=1}^N \mathbb{I}(x_i=j)\) for each outcome \(j \in \{1,\ldots,k\}\), where \(\mathbb{I}(\cdot)\) denotes the indicator function, and parameters \(\ak{\thetaN}=(\theta_{1,\ak{\thetaidx}},\ldots,\theta_{k,\ak{\thetaidx}})\in\mathbb{R}^k\) defining log-probability ratios \(\theta_{i,\ak{\thetaidx}}-\theta_{j,\ak{\thetaidx}} =\log [P_{\ak{\thetaN}}(X_1=i)/P_{\ak{\thetaN}}(X_1=j)]\). In addition to such standard models for discrete independent data, exponential models of FOES type commonly arise with dependent spatial data \citep{besag1974spatial} and network/relational data \citep{wasserman1994social, handcock2003assessing}. For a random graph or network with, say, \(n\) nodes \ak{$\{v_i\}_{i = 1}^n$}, consider \(N={n \choose 2}\) random edges where the \(i\)th edge is associated with a pair of nodes \(s_i \equiv \{v_{i_1},v_{i_2}\}\) and a binary variable \(X_i\in\{0,1\}\) indicating presence/absence of an edge among the node pair \(s_i\), \(i=1,\ldots,N\). Here the length \(N\) of the edge variable sequence \(\boldsymbol X = (X_1,\ldots,X_N)\) increases as a function of node number \(n\) and corresponding exponential models often incorporate graph topographical features derived from \(\boldsymbol X\). As an example, consider a graph model of exponential/FOES form prescribed by
\begin{equation}
\label{eq:mod2}
P_{\ak{\thetaN}}(\boldsymbol x) \propto
 \exp\left[\sum_{j=1}^3 \theta_{j,\ak{\thetaidx}} g_{j,N}(\boldsymbol x)\right], \quad\qquad \boldsymbol x=(x_1,\ldots,x_N)  \in\{0,1\}^N,
\end{equation}
\[
  g_{1,N}(\boldsymbol x) \equiv \sum_{i=1}^N  x_i, \qquad\quad g_{2,N}(\boldsymbol x) \equiv \sum_{1 \leq i<j \leq N,\atop s_i \cap s_j \neq \emptyset}\!\!\!   x_i x_j, \qquad
  g_{3,N}(\boldsymbol x) \equiv \sum_{1 \leq i<j<\ell \leq N, \atop s_i \cap s_j \neq \emptyset,s_i \cap s_\ell \neq \emptyset, \\ s_j \cap s_\ell \neq \emptyset } \!\!\!\!\!\!\!\!  \!\!\!\!\!\!\!\!   x_i x_j x_\ell,
\]
involving the numbers of edges, 2-stars and triangles among an outcome \(\boldsymbol x\) given by \(g_{1,N}(\boldsymbol x)\), \(g_{2,N}(\boldsymbol x)\) and \(g_{3,N}(\boldsymbol x)\), respectively, along with \(k=3\) real parameters \(\ak{\thetaN} \equiv (\theta_{1,\ak{\thetaidx}},\theta_{2,\ak{\thetaidx}},\theta_{3,\ak{\thetaidx}})\). For this network model \eqref{eq:mod2} in particular, as well as for more general models of form \eqref{eq:expo}, \citet{schweinberger2011instability} considered instability in such exponential models with sequences of fixed parameters \(\ak{\thetaN} = (\theta_1,\ldots,\theta_k)\in\mathbb{R}^k\), \(N \geq 1\), of fixed dimension \(k\).

For model sequences \(P_{\ak{\thetaN}}(\boldsymbol x)\equiv p_{N,\ak{\thetaN}}(\boldsymbol x)\) of the exponential type \eqref{eq:expo}, such as those in \eqref{eq:mod1}-\eqref{eq:mod2}, note that the dimension \(k\) of the parameter \(\ak{\thetaN}\in\Theta \subset \mathbb{R}^k\) necessarily remains the same for all sample sizes \(N \geq 1\) as the form of the natural parameter function \(\eta(\cdot)\) in \eqref{eq:expo} and the number of sufficient statistics \(\boldsymbol g_{N}(\boldsymbol x)\) do not depend on \(N\). Consequently, \(\ak{\thetaN}\) lies in a parameter space of fixed Euclidean dimension \(k\). However, this aspect need not be true for other types of FOES models considered in Sections \ref{rbm} - \ref{deep-learning}, where instead the numbers of parameters and sufficient statistics commonly increase with the sample size \(N\).

\hypertarget{rbm}{%
\subsection{Restricted Boltzmann Machines}\label{rbm}}

A restricted Boltzmann machine (RBM) is an undirected graphical model specified for discrete or continuous random variables, with binary variables being most common \citep[cf.~][]{smolensky1986information}. A RBM architecture has two layers, hidden (\(\mathcal{H}\)) and visible (\(\mathcal{V}\)), with conditional independence within each layer. Let \(\boldsymbol X = (X_1,\ldots,X_N)\) denote the \(N\) random variables for visibles with support \(\mathcal{X}^N\) and \(\boldsymbol H = (H_1,\ldots,H_{N_\mathcal{H}})\) denote the \(N_\mathcal{H}\) random variables for hiddens with support \(\mathcal{X}^{N_\mathcal{H}}\) where \(\mathcal{X} = \{-1,1\}\). For \ak{specified} parameters \(\ak{\thetaN}^{\mathcal{H}} \in \mathbb{R}^{N_\mathcal{H}}\), \(\ak{\thetaN}^{\mathcal{V}}\in \mathbb{R}^N\), and \(\ak{\thetaN}^{\mathcal{HV}}\) as a real matrix with dimension \(N_\mathcal{H} \times N\), the RBM model for \(\tilde{\boldsymbol X}=(\boldsymbol X,\boldsymbol H)\) has the joint probability mass function
\begin{equation}
\label{eq:RBM1}
\tilde{P}_{\ak{\thetaN}} (\tilde{\boldsymbol x}) = \exp\left[ (\ak{\thetaN}^{\mathcal{H}})^T \boldsymbol h + \boldsymbol (\ak{\thetaN}^{\mathcal{V}})^T \boldsymbol x + \boldsymbol h^T  \ak{\thetaN}^{\mathcal{HV}} \boldsymbol x - \psi(\ak{\thetaN})\right], \quad \tilde{\boldsymbol x} = (\boldsymbol x, \boldsymbol h) \in \{\pm 1\}^{N+N_\mathcal{H}}
\end{equation}
with \ak{normalizing function $\psi(\ak{\thetaN})$}.
Let \(\ak{\thetaN} = (\ak{\thetaN}^{\mathcal{H}}, \ak{\thetaN}^{\mathcal{V}}, \ak{\thetaN}^{\mathcal{HV}} ) \in \Theta_\ak{\thetaidx} \equiv \mathbb{R}^{q(N)}\), with \(q(N) = N + N_\mathcal{H} + N*N_\mathcal{H}\), denote \ak{a given} parameter vector for the RBM, as indexed by the number \(N\) of visible random variables (which may differ from the actual lengths of these parameter vectors). The probability mass function for the visible variables \(\boldsymbol X = (X_1, \dots, X_N)\) follows from marginalizing the joint specification to yield
\begin{equation}
\label{eq:RBM2}
P_{\ak{\thetaN}} (\boldsymbol x) = \sum\limits_{\boldsymbol h \in \{\pm 1\}^{N_{\mathcal{H}}}} \tilde{P}_{\ak{\thetaN}} (\boldsymbol x, \boldsymbol h), \qquad \boldsymbol x \in \{\pm 1\}^{N}\equiv \mathcal{X}^N.
\end{equation}
Here the baseline model \eqref{eq:RBM1} for hidden/visible variables is a linear exponential one in sufficient statistics \((\tilde{\boldsymbol X}, \boldsymbol X^T\boldsymbol H)\) using \(\tilde{\boldsymbol X}=(\boldsymbol X,\boldsymbol H)\) from \eqref{eq:RBM1}, but the form differs from the previous exponential models in \eqref{eq:expo} in that the lengths of parameters \(\ak{\thetaN}\) and statistics \((\tilde{\boldsymbol X}, \boldsymbol X^T\boldsymbol H)\) increase to incorporate more visible variables. That is, in contrast to \eqref{eq:expo}, the natural parameter function involved in the RBM model \eqref{eq:RBM1}, as the identity mapping of the parameters \(\ak{\thetaN}\in\mathbb{R}^{q(N)}\), naturally grows in dimension \(q(N)\to \infty\) to accommodate visible variables \(X_1, \dots, X_N\) of increasing data size \(N\to \infty\). Additionally, one may further arbitrarily choose the number \(N_\mathcal{H}\) of hidden variables \(\boldsymbol H\) in the joint RBM model \eqref{eq:RBM1} to define a marginal model \eqref{eq:RBM2} for the \(N\) visible variables \(\boldsymbol X\), and the number \(N_\mathcal{H}\) of hiddens may also potentially increase with \(N\). Because \(|\mathcal{X}| = 2\) and \(P_{\ak{\thetaN}}(\boldsymbol x) > 0\) for all \(\boldsymbol x \in \mathcal{X}^N\), the RBM specification \eqref{eq:RBM2} for visibles \(\boldsymbol\) corresponds to a FOES model, while the joint distribution \eqref{eq:RBM1} for \((\boldsymbol X, \boldsymbol H)\) is also a FOES model. As this example also indicates, any model formed by marginalizing a base FOES model class, such as the RBM joint specification \eqref{eq:RBM1}, is again a FOES model.

\ak{\subsection{Binary Deep Learning Models} \label{deep-learning}}

Consider two models with ``deep architecture'' that contain multiple hidden (or latent) layers in addition to a visible layer of data, namely a deep Boltzmann machine \citep{salakhutdinov2009deep} and a deep belief network \citep{hinton2006fast}. Let \(M\) denote the number of hidden layers included in the model and let \(N_{(H,1)}, \dots, N_{(H,M)}\) denote the numbers of hidden variables within each hidden layer. Then the random vector \(\tilde{\boldsymbol X} = \{H^{(1)}_1, \dots, H^{(1)}_{N_{(H,1)}}, \dots, H^{(M)}_1, \dots, H^{(M)}_{N_{(H,M)}}, \boldsymbol X\}\) collects both the hidden variables \(\{ H_{i}^{(j)} : i=1,\ldots, N_{(H,j)}, j=1,\ldots,M\}\) and visible variables \(\boldsymbol X =(X_1,\ldots,X_N)\) in a deep probabilistic model. Each variable outcome will again lie in \(\mathcal{X} = \{-1,1\}\).

\noindent \emph{Deep Boltzmann machine (DBM).} The DBM class of models maintains conditional independence within all layers in the model by stacking RBM models and only allowing conditional dependence between neighboring layers. The joint probability mass function for a DBM is
\[
\tilde{P}_{\ak{\thetaN}} ( \tilde{\boldsymbol x} ) = \exp\left[ \sum\limits_{i = 1}^M\boldsymbol \alpha^{(i)T} \boldsymbol h^{(i)} + \boldsymbol \beta^T \boldsymbol x + \boldsymbol h^{(1)T} \Gamma^{(0)} \boldsymbol x + \sum\limits_{i = 1}^{M - 1} \boldsymbol h^{(i)T} \Gamma^{(i)} \boldsymbol h^{(i + 1)} - \psi(\ak{\thetaN}) \right],
\]
for \(\tilde{\boldsymbol x} = (\boldsymbol h^{(1)}, \dots, \boldsymbol h^{(M)}, \boldsymbol x) \in \mathcal{X}^{N_{(H,1)} + \cdots + N_{(H,M)} +N}\) \ak{where again $\psi(\ak{\thetaN})$}
is the normalizing function for \(\ak{\thetaN} = (\boldsymbol \alpha^{(1)}, \dots, \boldsymbol \alpha^{(M)}, \boldsymbol \beta,\Gamma^{(0)}, \dots, \Gamma^{(M - 1)}) \in \Theta_\ak{\thetaidx} \subset \mathbb{R}^{q(N)}\), consisting of model parameters \(\boldsymbol \beta \in \mathbb{R}^N\), \(\boldsymbol \alpha^{(i)} \in \mathbb{R}^{N_{(H,i)}}\), \(i = 1, \dots, M\), along with a matrix \(\Gamma^{(0)}\) of dimension \(N_{(H,1)} \times N\), and matrices \(\Gamma^{(i)}\) of dimension
\(N_{(H,i)} \times N_{(H,i+1)}\) for \(i = 1, \dots, M-1\). The combined parameter vector \(\ak{\thetaN}\) has total length \(q(N)= N_{(H,1)}+\cdots N_{(H,M)} + N + N_{(H,1)}*N+N_{H,2}*H_{(H,1)}+\cdots +N_{(H,M)}*H_{(H,M)-1}\). The probability mass function for the visible random variables \(X_1, \dots, X_N\) follows from this joint specification as
\[
P_{\ak{\thetaN}} (\boldsymbol x) = \sum\limits_{(\boldsymbol h^{(1)}, \dots, \boldsymbol h^{(M)}) \in \mathcal{X}^{N_{(H,1)} + \cdots + N_{(H,M)}}} \tilde{P}_{\ak{\thetaN}} (\boldsymbol h^{(1)}, \dots, \boldsymbol h^{(M)}, \boldsymbol x) , \qquad \boldsymbol x \in \mathcal{X}^N.
\]
Again like the RBM case, the DBM model specification is an example of a FOES model.

\noindent \emph{Deep belief network (DBN).} A DBN resembles a DBM in that there are multiple layers of latent random variables stacked in a deep architecture with no conditional dependence between layers. The difference between the DBM and DBN models is that all but the last stacked layer \ak{(e.g. "layer $0$" corresponding to the visible random variables $X_1, \dots, X_N$)} in a DBN are Bayesian networks \citep[see][]{pearl985bayesian}, rather than RBMs. A Bayesian network \ak{is a class of probabilistic graphical models that define} conditional dependence to be directed, rather than undirected (as with the RBM). Thus for visibles \(X_1, \dots, X_N\) with support \(\mathcal{X}^N, \mid \mathcal{X} \mid < \infty\), a DBN is also a FOES model with \(q(N)\) the length of parameter vector is dependent on the dimension of the visibles because \(P_{\ak{\thetaN}}(\boldsymbol x)>0\) for all \(\boldsymbol x \in\mathcal{X}^N\). Commonly, as in logistic belief nets \citep{neal1992connectionist}, a ``weight'' parameter is placed on each interaction between visibles, \(X_1, \dots, X_N\), and the first layer of latent variables, \(H^{(1)}_1, \dots, H^{(1)}_{N_{(H,1)}}\), satisfying the definition of a FOES model.

\hypertarget{instability-results}{%
\section{Main Results on Model Instability}\label{instability-results}}

We now present a formal definition for instability of FOES models as well as a simple condition for identifying instability in a FOES model sequence. \ak{The basic intuition is that, if the smallest and largest probabilities from a given model are ``too different," then instability may arise, which manifests as extreme sensitivity in other probability behavior: large shifts in probability may be associated with only small differences in the input space for FOES models, and is related to a particular model placing almost all probability on a subset of potential outcomes. We present both non-asymptotic (Theorem \ref{thm:instab-elpr}) and asymptotic (Theorem \ref{thm:degenFOES}) results to this effect.}

\hypertarget{criterion}{%
\subsection{A Criterion for Instability}\label{criterion}}

To define a measure of instability in FOES models, it is useful to consider the behavior of data models \(P_{\ak{\thetaN}}\), again supported on a set \(\mathcal{X}^N\) of outcomes for \(\boldsymbol X\equiv \boldsymbol X_N =(X_1,\ldots,X_N)\) and for a prescribed configuration of parameters \(\ak{\thetaN}\), in connection to the sample size \(N\). A relevant quantity to this end is a log-ratio of extremal probabilities (LREP), defined as
\begin{align}
\label{eq:elpr}
 \REP (\ak{\thetaN})  =  \log \left[\frac{\max\limits_{  \boldsymbol x\in \mathcal{X}^N}P_{\ak{\thetaN}}( \boldsymbol x)}{\min\limits_{ \boldsymbol x \in \mathcal{X}^N}P_{\ak{\thetaN}}( \boldsymbol x)}\right],
\end{align}
based on maximum and minimal model probabilities. In what follows, the main idea is that instability, and other negative model features, can be associated with a FOES model formulation for \(N\) random variables where the LREP \eqref{eq:elpr} is overly large relative to the sample size \(N\). That is, a sequence of FOES probability models \(P_{\ak{\thetaN}}\) results in specifying \ak{a} distribution of observations \(\boldsymbol X=(X_1,\ldots,X_N)\) for each sample size \(N \geq 1\) \ak{(dependent on the prescribed values for $\ak{\thetaN}$, $N \geq 1$)} and instability will generally occur among these models whenever the corresponding LREP \eqref{eq:elpr} grows faster than \(N\). This leads to the following definition.

\BeginKnitrBlock{definition}[S-unstable FOES model]
\protect\hypertarget{def:instabFSFS}{}{\label{def:instabFSFS} \iffalse (S-unstable FOES model) \fi{} }A FOES model formulation \ak{$P_{\ak{\thetaN}}$} for \(\boldsymbol X_N=(X_1,\ldots,X_N)\), \ak{based on parameter specification $\ak{\thetaN}$, $N \geq 1$,} is \emph{Schweinberger-unstable} or \emph{S-unstable} if
\begin{equation}
\label{eq:Sun}
\lim \limits_{N \rightarrow \infty} \frac{1}{N} \REP(\ak{\thetaN}) \equiv \lim \limits_{N \rightarrow \infty} \frac{1}{N}\log \left[\frac{\max\limits_{  \boldsymbol x\in \mathcal{X}^N}P_{\ak{\thetaN}}( \boldsymbol x)}{\min\limits_{ \boldsymbol x \in \mathcal{X}^N}P_{\ak{\thetaN}}( \boldsymbol x)}\right] = \infty
\end{equation}
as the number of variables increases (\(N \rightarrow \infty\)).
\EndKnitrBlock{definition}

In other words, a model is S-unstable if \(\REP(\ak{\thetaN})/N\) is an unbounded sequence of sample size \(N\); namely, given any \(C > 0\), there exists an integer \(N_C > 0\) so that \(\REP(\ak{\thetaN})/N > C\) holds for all \(N \ge N_C\). A FOES model formulation may be termed S-stable if it fails to be S-unstable, i.e., if \(\sup_{N \geq 1}\REP(\ak{\thetaN})/N\) is bounded.

This definition of S-unstable is a generalization or reinterpretation of ``unstable'' used in \citet{schweinberger2011instability} by allowing possibly non-exponential family models (e.g., RBM and DBM models in Sections \ref{rbm}-\ref{deep-learning} as well as a potentially increasing number \(q(N)\) of parameters through the parameter sequence \(\ak{\thetaN}\in \mathbb{R}^{q(N)}\). While this definition differs in form and scope from the original, it does match that in \citet{schweinberger2011instability} for the special case of exponential models (cf.~Section \ref{discrete-exponential-family-models} considered there \ak{which translates to $\ak{\thetaN}=\boldsymbol{\theta} \in \mathbb{R}^k$, $N \geq 1$, in our model notation). Note the definition (\ref{eq:Sun}) of model instability depends intricately on the prescribed parameter sequence $\ak{\thetaN}$ or the particular parameter values indexing the distribution $P_{\ak{\thetaN}}$. This also agrees with the notion of unstable distributions from \citet{schweinberger2011instability}, as well as other characterizations of distributional stability in networks (cf. \citet{handcock2003assessing}), that model instability is crucially tied to the parameter configuration (i.e., where $\ak{\thetaN}$ lies in the parameter space) in addition to distributional form.  Hence, a given model family or form may may result in an unstable/stable model, depending on the parameters $\ak{\thetaN}$ chosen.} Section \ref{illustrations} provides several examples of unstable models as well as causes for model instability, where the latter may often be traced to issues in model form (i.e., data functions) and/or \ak{type of} parameterization. We next describe several potentially undesirable features associated with S-unstable FOES models.

\BeginKnitrBlock{remark}
{}\ak{Different notions of model stability/instability exist. For example, \citet{handcock2003assessing} refers to how small changes in parameters may dramatically change a probability mass function, while \citet{kaiser2007statistical} discusses how interpretation of model parameters and conditional expectations may break down in parts of the parameter space. However, these may be intuitively connected in the sense that the same (undesirable) parameter/model configurations may often be associated with varying conceptualizations of instability. While beyond the scope of the current paper, \citet{kaplan2019properties} investigates and numerically demonstrates some such connections for a class of RBM models. Section \ref{stable-models} also explains one relation between our notion of stable models and parameter interpretation in centered models \citep[][cf.]{kaiser2009exploring}. }
\EndKnitrBlock{remark}

\BeginKnitrBlock{remark}
{}In the definition \eqref{eq:Sun} of S-instability, we note that the numerical measure \(\REP(\ak{\thetaN})/N\) of model instability is invariant to \emph{independent replications} of data. \ak{Consequently, the definition of an
S-unstable model is unaffected by independent replication and all instability properties may be
characterized by those of one observation from the common FOES model. Remark A.2 of the Appendix provides more details.}
\EndKnitrBlock{remark}

\hypertarget{characterizations-and-consequences-of-instability}{%
\subsection{Characterizations and Consequences of Instability}\label{characterizations-and-consequences-of-instability}}

As a basic characteristic, S-unstable FOES model sequences have extremely sensitive probability structures. One aspect is that small changes in data configuration can lead to very large changes in probability. Consider, for example, the quantity given by
\[
\DN(\ak{\thetaN}) \equiv \max \left\{\log \frac{P_{\ak{\thetaN}}(\boldsymbol x)}{P_{\ak{\thetaN}}(\boldsymbol x^*)} : \boldsymbol x \text{ }\& \text{ } \boldsymbol x^* \in \mathcal{X}^N \text{ differ in exactly one component}\right\},
\]
which represents the biggest log-probability ratio for a one-component change \ak{(i.e., the data vector differs in exactly one component value)} in data outcomes in a FOES model with parameter \(\ak{\thetaN}\). \ak{This is the smallest difference possible between two data vectors, making it a good notion of a small distance in the input space and has been similarly used by \citet[][, Sec 2.2]{schweinberger2011instability}.} We then have the following result prescribing the behavior of \(\DN(\ak{\thetaN})\) for S-unstable FOES models.

\BeginKnitrBlock{theorem}
\protect\hypertarget{thm:instab-elpr}{}{\label{thm:instab-elpr} }Let \(P_{\ak{\thetaN}}\), with support \(\mathcal{X}^N\), \(N\geq 1\), be a sequence of FOES models.
\begin{enumerate}[(i)]
\item For any integer $N \geq 1$ and any given $C>0$, if $\REP(\ak{\thetaN})/N > C$ holds in (\ref{eq:elpr}), then
    $$ \DN(\ak{\thetaN}) > C$$
    follows.
\item Suppose the FOES model sequence is S-unstable. Then, \ak{for any $C>0$ there exists $N'$ such that for all $N>N'$}, there exist outcomes $\boldsymbol x,\boldsymbol x^*\in\mathcal{X}^N$, differing by one component, such that
    $$
    \frac{P_{\ak{\thetaN}}(\boldsymbol x)}{P_{\ak{\thetaN}}(\boldsymbol x^*)} > \exp[N C].
    $$
\end{enumerate}
\EndKnitrBlock{theorem}

\ak{Theorem \ref{thm:instab-elpr} aims to describe some model implications when a log-ratio of extreme model probabilities (\ref{eq:elpr}) is too large relative to the associated sample size $N$, which underlies the definition of a S-unstable FOES model (\ref{eq:Sun}). If so, Theorem \ref{thm:instab-elpr}(i) guarantees the FOES model must also exhibit correspondingly large changes in probability for small differences among some data configurations. Hence, as a further consequence in Theorem \ref{thm:instab-elpr}(ii), S-unstable models can never have universally bounded changes in probability among single component variations in data configurations. While not all one-component changes in data may produce massive changes in probability, unstable models must have some such data outcomes with this property. Hence,
unstable probability structures may exhibit extreme sensitivity through large peaks and troughs over the sample space.  Theorem \ref{thm:instab-elpr} supports and generalizes findings in \citet[Theorem 1]{schweinberger2011instability} for classes of discrete exponential models (i.e., the result there entails that, if $\REP(\thetaN)/N $ is unbounded as $N$ increases, then the nearest-point log-odds ratio  $\DN(\ak{\thetaN})$  is unbounded as well in these models).}

Additionally, S-unstable FOES model sequences are also connected to degenerate models, where \emph{degeneracy} involves assigning essentially all probability to modes within the sample space, which could potentially represent a small subset among the totality of outcomes. For perspective, note that differing sizes of the scaled log-ratio \(\REP(\ak{\thetaN})/N\) from \eqref{eq:elpr} induce a spectrum of levels of instability/stability and Theorem \ref{thm:instab-elpr} indicates increasing sensitivity of model probabilities as \eqref{eq:elpr} increases. Furthermore, as the instability measure grows and the log-ratio \(\REP(\ak{\thetaN})/N\) diverges, as in the definition \eqref{eq:Sun} of S-unstable models, then a FOES model sequence will become degenerate. Theorem \ref{thm:degenFOES} provides a formal statement of such degeneracy due to S-instability. For a given \(0 < \epsilon < 1\), define a \(\epsilon\)-modal set of outcomes as
\begin{equation}
\label{eq:mode}
\mathcal{M}_{\epsilon, \ak{\thetaN}} \equiv \left\{\boldsymbol x \in \mathcal{X}^N: \log P_{\ak{\thetaN}}(\boldsymbol x) > (1-\epsilon)\max\limits_{\boldsymbol y \in \mathcal{X}^N} \log  P_{\ak{\thetaN}}(\boldsymbol y) + \epsilon\min\limits_{\boldsymbol y \in \mathcal{X}^N} \log P_{\ak{\thetaN}}(\boldsymbol y) \right\}.
\end{equation}
\BeginKnitrBlock{theorem}
\protect\hypertarget{thm:degenFOES}{}{\label{thm:degenFOES} }For any arbitrarily small \(0 < \epsilon < 1\), an S-unstable FOES model sequence \(P_{\ak{\thetaN}}\), \(N \geq 1\), for \(\boldsymbol X_N=(X_1, \dots, X_N)\) satisfies
\[
P_{\ak{\thetaN}}\left( \boldsymbol X_N\in \mathcal{M}_{\epsilon, \ak{\thetaN}}\right) \rightarrow 1 \text{ as } N \rightarrow \infty.
\]
\EndKnitrBlock{theorem}

In other words, as the sample size grows in S-unstable FOES models, all probability tends to concentrate mass on an \(\epsilon\)-modal set, where \(\epsilon\) can be made arbitrarily small. Intuitively, the occurrence of such degeneracy can be explained by a type of ``reverse'' pigeonhole principle for unstable FOES models: if all outcomes should receive positive probability but the maximal probability far exceeds the minimal one in the model, then little probability remains for distribution among remaining model outcomes (i.e., if nearly all available pigeons are stuffed into one hole, the remaining pigeonholes must have few occupants). Degeneracy in unstable models can pose dangers in data modeling as well, particularly when a mode set represents a narrow collection of outcomes among those realistically possible for adequately describing data. In which case, model outcomes may fail to look like data of interest.

Connected to degeneracy, S-unstable FOES models may also exhibit additional kinds of extreme and undesirable sensitivity in probabilities if model parameters \(\ak{\thetaN}\) can further be ``dialed'' between positive and negative values. That is, some FOES models naturally involve parameter spaces covering a positive-negative spectrum of parameter possibilities, where the signs of parameters provide a standard device for increasing or decreasing probabilities of outcomes in the model formulation. In fact, for many models, the switch of a parameter sign serves to produce reciprocal probabilities, as outlined in the following model assumption about parameter sign reversal (PSR).

\textbf{Model Condition PSR} \emph{(Reciprocal Probabilities from Parameter Sign Reversal)}: Let \(P_{\ak{\thetaN}}\), with support \(\mathcal{X}^N\), \(N\geq 1\), represent a sequence of FOES models. For each \(N \geq 1\) and any outcome \(\boldsymbol x \in \mathcal{X}^N\), suppose it holds that
\[
P_{\ak{\thetaN}}(\boldsymbol x)  \cdot P_{-\ak{\thetaN}}(\boldsymbol x) \;=\;   \max\limits_{\boldsymbol y \in \mathcal{X}^N}P_{ \ak{\thetaN}}(\boldsymbol y)\cdot \min\limits_{\boldsymbol y \in \mathcal{X}^N}P_{-\ak{\thetaN}}(\boldsymbol y),
\]
where \(\max_{\boldsymbol y \in \mathcal{X}^N}P_{\ak{\thetaN}}(\boldsymbol y)\) and \(\min_{\boldsymbol y \in \mathcal{X}^N}P_{-\ak{\thetaN}}(\boldsymbol y)\) denote the maximum and minimum probabilities under parameters \(\ak{\thetaN}\) and \(-\ak{\thetaN}\), respectively.

The above model condition incorporates many standard parameterizations and follows, for instance, whenever \(P_{\ak{\thetaN}}(\boldsymbol x)/P_{\ak{\thetaN}}(\boldsymbol y) = [P_{-\ak{\thetaN}}(\boldsymbol x)/P_{-\ak{\thetaN}}(\boldsymbol y)]^{-1}\) holds for outcomes \(\boldsymbol x, \boldsymbol y \in\mathcal{X}^N\) in a FOES model. For instance, this latter condition is fulfilled for all linear exponential families from Section \ref{discrete-exponential-family-models} (e.g., \eqref{eq:mod1}-\eqref{eq:mod2}) as well as all network models from Sections \ref{rbm}-\ref{deep-learning} (e.g., \eqref{eq:RBM1}-\eqref{eq:RBM2}). When parameters can be tuned in sign with effects prescribed in the model condition PSR, unstable FOES models will exhibit further probability sensitivities, as outlined in the following extension of Theorem \ref{thm:degenFOES}.

\BeginKnitrBlock{corollary}
\protect\hypertarget{cor:sign}{}{\label{cor:sign} }Let \(P_{\ak{\thetaN}}\), with support \(\mathcal{X}^N\), \(N\geq 1\), be a sequence of FOES models satisfying model condition PSR. If the models \(P_{\ak{\thetaN}}\) are additionally S-unstable, then
\begin{enumerate}[(i)]
\item the models  $P_{-\ak{\thetaN}}$  defined by $-\ak{\thetaN}$  are also S-unstable;
\item and for the complement $\mathcal{M}_{\epsilon, \ak{\thetaN}}^c \equiv \mathcal{X}^N \setminus \mathcal{M}_{\epsilon, \ak{\thetaN}}$ of any mode-set $\mathcal{M}_{\epsilon, \ak{\thetaN}}$ under $\ak{\thetaN}$ from (\ref{eq:mode}), with $0<\epsilon<1$, it holds under $-\ak{\thetaN}$ that
    $$ 
    P_{-\ak{\thetaN}}( \boldsymbol X_N \in \mathcal{M}_{\epsilon, \ak{\thetaN}}^c) \rightarrow 1\qquad \text{as } N\to \infty,
    $$
    while, by Theorem \ref{thm:degenFOES}, $P_{\ak{\thetaN}}( \boldsymbol X_N \in \mathcal{M}_{\epsilon, \ak{\thetaN}} ) \rightarrow 1$ holds for $\boldsymbol X_N = (X_1,\ldots,X_N)$ under $\ak{\thetaN}$.
\end{enumerate}
\EndKnitrBlock{corollary}

For unstable models, Corollary \ref{cor:sign} shows that shifts in parameters around zero (i.e., from \(\ak{\thetaN}\) to \(-\ak{\thetaN}\)) can induce extreme changes in probability among subsets of the sample space, as another manifestation of instability and hyper-sensitivity in probability structure. For one-parameter exponential families, involving a fixed real-valued linear parameter \(\ak{\thetaN} = \theta \in \mathbb{R}\) and sufficient statistic \(\boldsymbol g_N(\boldsymbol x)\in \mathbb{R}\) in \eqref{eq:expo}, \citet[Theorem 3]{schweinberger2011instability} proved a result similar in spirit, though based on a characterization there in terms of maximum \(U_N \equiv \max_{\boldsymbol x\in\mathcal{X}^N}g_N(\boldsymbol x)\) and minimal \(L_N \equiv \min_{\boldsymbol x\in\mathcal{X}^N}g_N(\boldsymbol x)\) values of the sufficient statistic. For this case in particular, mode sets have specific, and essentially complementary, forms over positive and negative parameters, namely, \(\mathcal{M}_{\epsilon, \ak{\thetaN}} = \{\boldsymbol x \in\mathcal{X}^N: g_N(\boldsymbol x) > (1-\epsilon) U_N + \epsilon L_N \}\) and \(\mathcal{M}_{\epsilon, -\ak{\thetaN}} = \{\boldsymbol x \in\mathcal{X}^N: g_N(\boldsymbol x) < \epsilon U_N + (1-\epsilon) L_N \}\) for any \(\ak{\thetaN}>0\), and
\citet[Theorem 3]{schweinberger2011instability} showed each mode set collects all mass, under positive and negative parameters, respectively, with unstable models of this exponential type. However, for all unstable FOES models, Corollary \ref{cor:sign} generalizes the same principle that unstable models can push all probability to different, and in fact disjoint, parts of the sample space, depending on how parameters fall with respect to zero. This feature can numerically complicate likelihood manipulations, such as maximization or MCMC-based Bayes posterior \ak{simulation}, as further discussed in Section \ref{implications}.

\BeginKnitrBlock{remark}
{}Under the model condition PSR, Corollary \ref{cor:sign} can also be extended to
cases where parameter components \(\ak{\thetaN} =(\boldsymbol \theta_{1,\ak{\thetaidx}}, \boldsymbol \theta_{2, \ak{\thetaidx}})\) (say) are not all changed in sign (e.g., \(-\ak{\thetaN}\)) but, more generally, are instead altered to another parameter configuration \(\ak{\thetaN}^A = ( \boldsymbol \theta_{1,\ak{\thetaidx}}^A, \boldsymbol \theta_{2,\ak{\thetaidx}}^A )\) involving a switch in sign only among some dominating model parameters \(\boldsymbol \theta_{2,\ak{\thetaidx}}^A=- \boldsymbol \theta_{2,\ak{\thetaidx}}\) with remaining parameters \(\boldsymbol \theta_{1,\ak{\thetaidx}}^A\) being arbitrarily chosen. If a change sign occurs among parameters (\(\pm \boldsymbol \theta_{2,\ak{\thetaidx}}\)) which dominate the probability structure of the model, then the results of Corollary \ref{cor:sign} can still hold with \(\ak{\thetaN}^A\) replacing \(-\ak{\thetaN}\); as an example of one sufficient condition, if \(\lim_{N\to \infty} \max_{\boldsymbol x \in \mathcal{X}^N} |G_N(\boldsymbol x, -\ak{\thetaN}) - G_N(\boldsymbol x, \ak{\thetaN}^A)|=0\) holds in addition to Corollary \ref{cor:sign} assumptions, where
\[
G_N(\boldsymbol x, \boldsymbol \theta) =   \frac{\log P_{\boldsymbol \theta}(\boldsymbol x) - \min_{\boldsymbol y \in\mathcal{X}^N}\log P_{\boldsymbol \theta}(\boldsymbol y) }{\max_{\boldsymbol y \in\mathcal{X}^N}\log P_{\boldsymbol \theta}(\boldsymbol y)  - \min_{\boldsymbol y \in\mathcal{X}^N}\log P_{\boldsymbol \theta}(\boldsymbol y) }, \quad \boldsymbol x \in\mathcal{X}^N,
\]
represents a standardized form of \(\boldsymbol \theta\)-model probabilities, then the results of Corollary \ref{cor:sign} apply to \(\ak{\thetaN}^A\) in addition to \(-\ak{\thetaN}\). As a consequence, an unstable model under \(\ak{\thetaN}\) can then imply that many more unstable models exist over a broader spectrum of possibilities for variations \(\ak{\thetaN}^A\) of \(\ak{\thetaN}\), which involves some amount of sign change among components of \(\ak{\thetaN}\).
\EndKnitrBlock{remark}

\hypertarget{illustrations}{%
\section{Illustrations}\label{illustrations}}

Model instability can depend intricately on how functions of parameters and data \(\boldsymbol X_N=(X_1,\ldots,X_n)\) are combined in the formulation of the model probabilities, though some general causes may be identified. As one issue, a broad parameter space (or wide interpretation of this space) may admit some parameters as technically valid that have an undue and often undesirable impact on the model structure for a prescribed data size \(N\). In this case, both the size and dimension of model parameters can be problematic and induce instability. In combination to this last point, further causes of instability may also be traced to the magnitude of statistics in the model. Potentially massive, and thereby unstable, statistics were the primary focus of instability studies of \citet{schweinberger2011instability} for certain discrete exponential models having parameters/statistics of fixed dimension. However, as shown in the following, bounded statistics may still lead to instability if the parameter dimension is high. We next \ak{provide} some examples to illustrate S-instability in FOES models, which also suggest some potential strategies for preventing unstable models \ak{(e.g., for many model types in the following, a control of the combined magnitudes of parameters can ensure S-stability, which is also discussed further in Section \ref{stable-models})}.

\hypertarget{equi-probability-models}{%
\subsection{Equi-probability Models}\label{equi-probability-models}}

As a baseline for comparisons, consider a simplistic model for \(\boldsymbol X_N=(X_1,\ldots,X_N)\) with uniform probabilities over the sample space, say \(P_{\ak{\thetaN}}(\boldsymbol x)= |\mathcal{X}|^{-N}\), \(\boldsymbol x \in \mathcal{X}^N\), where each random variable has \(|\mathcal{X}| \geq 1\) outcomes. In contrast to instability, model probabilities here are completely insensitive to changes in data outcomes across the sample space, and the associated log-ratio of extreme probabilities \eqref{eq:elpr} is
\[
\frac{1}{N} \REP(\ak{\thetaN})=0\quad \;(\text{uniform probability model}),
\]
which is as small as possible. In fact, a LREP value of zero can only occur for a FOES model having uniform probabilities, and such equi-probability models are always S-stable.

\hypertarget{one-param-exp}{%
\subsection{One-parameter Exponential Models}\label{one-param-exp}}

A fundamental model considered in the instability work of \citet{schweinberger2011instability} involves a one-parameter exponential model corresponding to \eqref{eq:expo} with a real-valued parameter, say \(\ak{\thetaN} = \eta(\ak{\thetaN})\in \mathbb{R}\), and sufficient statistic \(g_N(\boldsymbol x)\in \mathbb{R}\).
For such models, upon scaling by sample size \(N\), the log-ratio of extreme probabilities in \eqref{eq:elpr} for assessing instability becomes
\begin{equation}
\label{eq:UL}
\frac{1}{N}\REP(\ak{\thetaN} ) \equiv   |\ak{\thetaN}| \frac{(U_N-L_N)}{N} \;\quad \text{(one-parameter exponential model)},
\end{equation}
where \(U_N \equiv \max_{\boldsymbol x\in\mathcal{X}^N}g_N(\boldsymbol x)\) and \(L_N \equiv \min_{\boldsymbol x\in\mathcal{X}^N}g_N(\boldsymbol x)\) denote the maximal and minimal values of the single sufficient statistic. In this case, an S-unstable model results, by definition \eqref{eq:Sun}, whenever \(\lim_{N\to \infty} |\ak{\thetaN}| (U_N-L_N)/N= \infty\) holds or, in other words, if the combined magnitudes of parameter \(|\ak{\thetaN}|\) and maximal difference \(U_N-L_N\) in statistic values are overwhelmingly large relative to the sample size \(N\). If we further assume that \(\ak{\thetaN} =\theta\in\mathbb{R}\setminus \{0\}\) is a fixed (non-zero) parameter for all \(N \geq 1\), as considered in \citet{schweinberger2011instability}, then an S-unstable model results solely if the sufficient statistic admits a value \(U_N-L_N\) too large relative to number \(N\) of observations, i.e., if \((U_N-L_N)/N\to \infty\) as \(N\to \infty\). The latter aspect reflects the definition of \citet{schweinberger2011instability}, for this setting, that a real-valued \emph{statistic} \(g_N(\boldsymbol x)\) may be classified as \emph{unstable} when \(\lim_{N\to \infty}|(U_N-L_N)/N=\infty\) holds and as \emph{stable} otherwise (e.g., if \(\sup_{N \geq 1}(U_N-L_N)/N<\infty\)).

For illustration, consider the iid Bernoulli model \eqref{eq:mod1} for \(\boldsymbol X_N=(X_1,\ldots,X_N)\) with log-odds ratio parameter \(\ak{\thetaN} = \log[ P_{\ak{\thetaN}}(X_1=1)/ P_{\ak{\thetaN}}(X_1=0)]\in\mathbb{R}\). \ak{Remark 3.2} (Section \ref{criterion}) then gives the model instability measure \eqref{eq:Sun} directly as
\[
\frac{1}{N}\REP(\ak{\thetaN} ) = |\ak{\thetaN}|\quad\; \text{(iid Bernoulli model)},
\]
so that an unstable (or stable) model results for a divergent (or bounded) parameter sequence \(|\ak{\thetaN}|\).
\ak{In this case}, \citet{schweinberger2011instability} has also noted that the statistic is stable (i.e., bounded \((U_N-L_N)/N=1\)) and the Bernoulli model is as well when, in particular, \(\ak{\thetaN}=\theta \in\mathbb{R}\) is fixed for \(N \geq 1\).

Alternatively, considering a random graph with \(N={n \choose 2}\) edges among \(n\) nodes, the exponential graph model from \eqref{eq:mod2}, when based purely on the number of \(g_{2,N}(\boldsymbol x)\) of 2-stars or solely the number \(g_{3,N}(\boldsymbol x)\) of triangles, \(\boldsymbol x\in\{0,1\}^N\), has a measure of instability from \eqref{eq:Sun} as
\[
\frac{1}{N}\REP(\ak{\thetaN} )  = \left\{ \begin{array}{lcl} |\ak{\thetaN}| (n-2) && \text{(2-star graph model)}\\
|\ak{\thetaN}|(n-2)/3 &&\text{(triangle graph model)},\end{array}\right.
\]
by using the (one-parameter exponential) LREP formula \eqref{eq:UL} with statistic maximums \(U_N= N(n-2)\) for 2-stars or \(U_N= N(n-2)/3\) for triangles and with minimums \(L_N=0\) in both cases. Because the variable number \(N\to \infty\) as the node number \(n\to \infty\), both counts of 2-stars and triangles are unstable statistics in the sense of \citet{schweinberger2011instability} (i.e., \(\lim_{N\to \infty} (U_N-L_N)/N=\infty\)). Furthermore, both types of graph models are always S-unstable for all possible of parameter sequences \(\ak{\thetaN} \in\mathbb{R}\) that are bounded away from zero (i.e., \(\lim_{N\to \infty}\REP(\ak{\thetaN} )/N=\infty\) then holds, including the fixed parameter case \(\ak{\thetaN}=\theta\in\mathbb{R}\setminus \{0\}\) from \citet{schweinberger2011instability}).

\hypertarget{fixed-dim-exp}{%
\subsection{Fixed-dimensional Linear Exponential Models}\label{fixed-dim-exp}}

As a generalization of the one-parameter exponential case, we next consider linear exponential families \eqref{eq:expo} with \(k\) parameters \(\ak{\thetaN} = (\theta_{1,\ak{\thetaidx}},\ldots,\theta_{k,\ak{\thetaidx}})^\prime\) and \(k\) sufficient statistics \(\boldsymbol g_N(\boldsymbol x) = (g_{1,N}(\boldsymbol x),\ldots, g_{k,N}(\boldsymbol x))^\prime\). Here the dimension \(k\) of model parameters/statistics is fixed, and we next prescribe a condition helpful to avoiding instability in such models. For this, define \(U_{i,N}=\max_{\boldsymbol x \in\mathcal{X}^N} g_{i,N}(\boldsymbol x)\) and \(L_{i,N}=\min_{\boldsymbol x \in\mathcal{X}^N} g_{i,N}(\boldsymbol x)\) as the maximal and minimal values of the \(i\)th statistic, \(i=1,\ldots,k\), based on observations \(\boldsymbol X_N=(X_1,\ldots,X_N)\).

\BeginKnitrBlock{proposition}
\protect\hypertarget{prp:prop1}{}{\label{prp:prop1} }Let \(P_{\ak{\thetaN}}\), \(N \geq 1\), denote linear exponential models (\ref{eq:expo}) with parameters \(\ak{\thetaN} = (\theta_{1,\ak{\thetaidx}},\ldots,\theta_{k,\ak{\thetaidx}})^\prime \in \mathbb{R}^k\) and statistics \(\boldsymbol g_N(\boldsymbol x) = (g_{1,N}(\boldsymbol x),\ldots, g_{k,N}(\boldsymbol x))^\prime \in \mathbb{R}^k\), for fixed \(k \geq 1\). Then, the models \(P_{\ak{\thetaN}}\) are S-stable if
\begin{equation}
\label{eq:prop1}
\sup_{N \geq 1}\frac{1}{N} \max_{1 \leq i \leq k }|\theta_{i,\ak{\thetaidx}}|(U_{i,N}-L_{i,N})<\infty
\end{equation}
holds, i.e., if \(\max_{1 \leq i \leq k } |\theta_{i,\ak{\thetaidx}}|(U_{i,N}-L_{i,N})/N\) is bounded sequence of sample size \(N\).
\EndKnitrBlock{proposition}

\BeginKnitrBlock{remark}
{}In the one-parameter exponential case \(k=1\), recall the exponential model is stable/unstable depending on whether \(\REP(\ak{\thetaN})/N = |\theta_{1,\ak{\thetaidx}}|(U_{1,N}-L_{1,N})/N \equiv | \ak{\thetaN}|(U_{N}-L_{N})/N\) in \eqref{eq:UL} is convergent/divergent. Hence, for \(k=1\), the condition \eqref{eq:prop1} of Proposition \ref{prp:prop1} captures the same notion of S-stability based on \eqref{eq:UL}.
\EndKnitrBlock{remark}

Proposition \ref{prp:prop1} provides a sufficient condition for the stability of linear exponential models with fixed parameter dimension \(k\geq 1\), whereby an S-stable model is guaranteed if the \ak{product of magnitudes} of each combination of parameter \(\theta_{i,\ak{\thetaidx}}\) and sufficient statistic value \((U_{i,N}-L_{i,N})\) is bounded by the sample size \(N\), \(i=1,\ldots,k\). This supports the findings of \citet{schweinberger2011instability}, who showed degeneracy follows in such models under one type of violation of the condition \eqref{eq:prop1} in Proposition \ref{prp:prop1} (namely, involving \(k>1\) non-zero parameters with \(k-1\) statistics being \(O(N)\) bounded while one statistic diverges in maximal size faster than the number \(N\) of observations). To further illustrate the result in Proposition \ref{prp:prop1}, consider the multinomial distribution \eqref{eq:mod11} for \(\boldsymbol X_N=(X_1,\ldots,X_N)\) having \(k\geq 2\) categories \(\{1,\ldots,k\}\) and \(k\) parameters \(\ak{\thetaN} = (\theta_{1,\ak{\thetaidx}},\ldots,\theta_{k,\ak{\thetaidx}})^\prime\). The variables are iid under this model so that \ak{Remark 3.2} (Section \ref{criterion}) yields the corresponding \(N\)-scaled log-ratio of extreme probabilities \eqref{eq:elpr} as
\begin{align*}
  \frac{1}{N}\REP(\ak{\thetaN}) &= \frac{\max_{1 \leq i \leq k} P_{\ak{\thetaN}}(X_1=i)}{\min_{1 \leq i \leq k} P_{\ak{\thetaN}}(X_1=i)}\\
  &= \max_{1 \leq i \leq k} \theta_{i,\ak{\thetaidx}} - \min_{1 \leq i \leq k} \theta_{i,\ak{\thetaidx}} \qquad \text{(iid multinomial model)}.
\end{align*}
Hence, a multinomial model sequence is unstable (or stable) depending on whether (or not) the maximal parameter difference \(\max_{1 \leq i \leq k} \theta_{i,\ak{\thetaidx}} - \min_{1 \leq i \leq k} \theta_{i,\ak{\thetaidx}}\) diverges. Furthermore, using that each of the \(k\) sufficient (count) statistics from the multinomial model \eqref{eq:mod11} satisfies \((U_{i,N}-L_{i,N})/N=1\), we see that \eqref{eq:prop1} of Proposition \ref{prp:prop1} becomes purely a parameter condition, \(\sup_{N \geq 1}\max_{1 \leq i \leq k } |\theta_{i,\ak{\thetaidx}}| <\infty\), for ensuring that \(\REP(\ak{\thetaN})/N =\max_{1 \leq i \leq k} \theta_{i,\ak{\thetaidx}} - \min_{1 \leq i \leq k} \theta_{i,\ak{\thetaidx}}\) is bounded and stability follows for the multinomial distribution. Additionally, a stable multinomial sequence (i.e., bounded \(\REP(\ak{\thetaN})/N\)) turns out to be nearly equivalent to \eqref{eq:prop1} (e.g., these are the same if the smallest parameter \(\min_{1 \leq i \leq k } |\theta_{i,\ak{\thetaidx}}|\) remains bounded).

When the condition \eqref{eq:prop1} of Proposition \ref{prp:prop1} is violated, this aspect suggests a potentially unstable model that may be investigated more closely. For example, consider the exponential graph model from \eqref{eq:mod2} involving counts of edges, 2-stars and triangles with fixed parameters \(\thetaN = (\theta_{1},\theta_2,\theta_3)^\prime \in \mathbb{R}^3\) for \(N\geq 1\). If either the 2-star parameter \(\theta_2 \neq 0\) or triangle parameter \(\theta_3 \neq 0\) is non-zero, then \(\max_{1 \leq i \leq 3 } |\theta_{i}|(U_{i,N}-L_{i,N})/N \propto (n-2)\to \infty\) holds in \eqref{eq:prop1} by \((U_{2,N}-L_{2,N})/N = 3 (U_{3,N}-L_{3,N})/N=(n-2)\) for 2-star and triangle statistics (\(i=2,3\)), so that Proposition \ref{prp:prop1} hints that an unstable model may result when \(|\theta_2| + |\theta_3| \neq 0\). Relatedly, a result from \citet[Result 3]{schweinberger2011instability} states that this model is unstable for all fixed parameters excluding cases \(\theta_2 =\theta_3=0\) or \(\theta_2 = - \theta_3/3\)\ak{. However, more in line} with the instability suggested by Proposition \ref{prp:prop1}\ak{, this model is unstable} whenever \(|\theta_2| + |\theta_3| \neq 0\) (i.e., excluding \(\theta_2 =\theta_3=0\))\ak{; Remark A.1 in the Appendix gives a proof.} That is, instability holds even under \ak{the} \(\theta_2 = - \theta_3/3\) case potentially allowed by Schweinberger's \citeyearpar{schweinberger2011instability} results. \ak{Thus, graph models of the form (\ref{eq:mod2}), with 2-stars and triangles are always S-unstable for all possible parameter sequences $\thetaN \in \mathbb{R}^3$ with $(\theta_2, \theta_3)$ bounded away from $\boldsymbol 0 \in \mathbb{R}^2$.}

\hypertarget{latent-variable-models-of-increasing-parameter-dimension}{%
\subsection{Latent Variable Models of Increasing Parameter Dimension}\label{latent-variable-models-of-increasing-parameter-dimension}}

We next consider instability of discrete data models based on exponential formulations involving hidden, or latent, variables, such as those probabilistic graphical models described in Sections \ref{rbm}-\ref{deep-learning}. We will focus on restricted Boltzmann machine (RBM) models (Section \ref{rbm}, having one layer of latent variables for simplicity, though the same instability concepts may be extended to other deep learning models (Section \ref{deep-learning}. For \(N\) visible variables \(\boldsymbol X \equiv \boldsymbol X_N = (X_1,\ldots,X_n)\) as data, each observation \(X_i\in\{\pm 1\}\) being binary, the RBM-based model \eqref{eq:RBM2} for \(\boldsymbol X\) is again of FOES-type, though not an exponential model. However, the distribution of visible variables is induced by an underlying joint exponential model \eqref{eq:RBM1} for both visible and latent variables \((\boldsymbol X, \boldsymbol H)\), where \(\boldsymbol H = (H_1,\ldots,H_{N_{\mathcal{H}}})\) denotes a vector of \(N_{\mathcal{H}}\) hidden variables (similarly binary). The joint model is of linear exponential form involving \(q(N)\equiv N + N_{\mathcal{H}} + N*N_{\mathcal{H}}\) sufficient statistics given by \((\boldsymbol X, \boldsymbol H, \boldsymbol X^T\boldsymbol H)\) and parameters \(\ak{\thetaN} = (\ak{\thetaN}^{\mathcal{V}}, \ak{\thetaN}^{\mathcal{H}}, \ak{\thetaN}^{\mathcal{VH}} ) \in\mathbb{R}^{q(N)}\) corresponding to the \(N\) visible variables \(\boldsymbol X\) (i.e., \(\ak{\thetaN}^{\mathcal{V}}\in\mathbb{R}^N\)), the \(N_{\mathcal{H}}\) hidden variables \(\boldsymbol H\) (i.e., \(\ak{\thetaN}^{\mathcal{H}}\in\mathbb{R}^{N_{\mathcal{H}}}\)), and the \(N *N_{\mathcal{H}}\) cross-product variables \(\boldsymbol X^T\boldsymbol H\) (i.e., \(\ak{\thetaN}^{\mathcal{VH}}\in\mathbb{R}^{N*N_{\mathcal{H}}}\)). However, unlike some previous exponential models considered in Sections \ref{one-param-exp}-\ref{fixed-dim-exp} (cf.~Proposition \ref{prp:prop1}, note that the RBM formulation always associates parameters with \emph{bounded} statistics (i.e., the components of \((\boldsymbol X, \boldsymbol H, \boldsymbol X^T\boldsymbol H)\)) so that model instability cannot arise here due to the magnitude of sufficient statistics exceeding the sample size \(N\). Instead, RBM instability may be linked solely to parameter configuration and the fact that the number \(q(N) \geq N\) of parameters necessarily increases with the number \(N\) of observations \(\boldsymbol X\), in contrast to previous exponential cases of fixed parameter dimension.

To highlight the instability issues for the RBM model, consider a simple model for \(N\) visibles \(\boldsymbol X\) with no hidden variables (\(N_{\mathcal{H}}=0\)), for which model statements \eqref{eq:RBM1}-\eqref{eq:RBM2} coincide. An independence model then results for variables in \(\boldsymbol X\), which has \(q(N)=N\) parameters \(\ak{\thetaN}^{\mathcal{V}} = (\theta_{1,\ak{\thetaidx}}^{\mathcal{V}}, \ldots, \theta_{N,\ak{\thetaidx}}^{\mathcal{V}}) \in \mathbb{R}^{N}\), and the measure of model instability becomes
\[
\frac{1}{N}\REP(\ak{\thetaN}) = \frac{2}{N} \sum_{i=1}^N|\theta_{i,\ak{\thetaidx}}^{\mathcal{V}}| \quad \text{(RBM model, no hiddens)}.
\]
Hence, this model sequence for \(\boldsymbol X\) will be S-unstable if the aggregation of absolute parameters grows faster than the number \(N\) of parameters/visible variables. Consequently, even for \ak{the} simplest RBM model involving independence, preventing instability requires careful choice of parameters, particularly with regard to how a parameter configuration differs from zero. For more general RBM models, the number \(N_{\mathcal{H}}\) of hidden variables \(\boldsymbol H\) can also be chosen arbitrarily (i.e., as some function \(N_{\mathcal{H}}\equiv N_{N,\mathcal{H}}\) of \(N\)), which can substantially inflate the number \(q(N)\) of model parameters and further impact model instability through accumulated parameters. To better understand the effects of instability in the RBM structure, Proposition \ref{prp:prop2} next frames the general behavior of extreme probabilities in the joint RBM model \eqref{eq:RBM1} for \((\boldsymbol X, \boldsymbol H)\) and the implied RBM data model \eqref{eq:RBM2} for \(\boldsymbol X\) alone. Specifically, critical measures of instability may be closely connected in both models through tight bounds on their respective LREP values \eqref{eq:elpr}. As a result, Proposition \ref{prp:prop2} shows how an unstable distribution for observations \(\boldsymbol X\) may be traced to sources of instability in the original joint distribution for \((\boldsymbol X,\boldsymbol H)\). This also suggests a device for avoiding instability, as provided next.

To state the result, let \(\REP_{\boldsymbol X}(\ak{\thetaN}) \equiv \REP(\ak{\thetaN})\) denote the LREP value \eqref{eq:elpr} from the marginal distribution \(P_{\ak{\thetaN}}\) of visibles \(\boldsymbol X\) in \eqref{eq:RBM2} and write the LREP for the joint distribution \(\tilde{P}_{\ak{\thetaN}}\) of \((\boldsymbol X, \boldsymbol H)\) from \eqref{eq:RBM1} as
\begin{align*}
\REP_{(\boldsymbol X, \boldsymbol H)}(\ak{\thetaN}) &=\log\left[  \frac{\max_{ (\boldsymbol x, \boldsymbol h) \in \{\pm 1\}^{N+N_{\mathcal{H}}}}
\tilde{P}_{\ak{\thetaN}} (\boldsymbol x, \boldsymbol h)}{\min_{ (\boldsymbol x, \boldsymbol h) \in \{\pm 1\}^{N+N_{\mathcal{H}}}`}
\tilde{P}_{\ak{\thetaN}} (\boldsymbol x, \boldsymbol h)} \right] \qquad \text{(joint RBM model)}    \\    
&= \left(\max_{ (\boldsymbol x, \boldsymbol h) \in \{\pm 1\}^{N+N_{\mathcal{H}}}} f_{\ak{\thetaN}} (\boldsymbol x, \boldsymbol h)-\min_{ (\boldsymbol x, \boldsymbol h) \in \{\pm 1\}^{N+N_{\mathcal{H}}}} f_{\ak{\thetaN}} (\boldsymbol x, \boldsymbol h)\right),
\end{align*}
written as a function
\begin{equation}
\label{eq:f}
f_{\ak{\thetaN}} (\boldsymbol x, \boldsymbol h)  \equiv  \sum_{i=1}^N x_i  \theta_{i,\ak{\thetaidx}}^{\mathcal{V}} +  \sum_{j=1}^{N_\mathcal{H}} h_j \theta_{j,\ak{\thetaidx}}^{\mathcal{H}} + \sum_{i=1}^N \sum_{j=1}^{N_\mathcal{H}} x_i h_j  \theta_{ij,\ak{\thetaidx}}^{\mathcal{VH}}  
\end{equation}
of outcomes \(\boldsymbol x =(x_1,\ldots,x_N) \in\{\pm1\}^{N}\), \(\boldsymbol h =(h_1,\ldots,h_{N_{\mathcal{H}}}) \in\{\pm1\}^{N_{\mathcal{H}}}\) and parameters \(\ak{\thetaN} \equiv (\ak{\thetaN}^{\mathcal{V}}, \ak{\thetaN}^{\mathcal{H}}, \ak{\thetaN}^{\mathcal{VH}})\), with \(\theta_{i,\ak{\thetaidx}}^{\mathcal{V}}\), \(\theta_{j,\ak{\thetaidx}}^{\mathcal{H}}\) and \(\theta_{ij,\ak{\thetaidx}}^{\mathcal{VH}}\) denoting respective parameter components, \(1 \leq i \leq N\), \(1 \leq j \leq N_{\mathcal{H}}\). Due to the marginalization steps in defining the distribution \eqref{eq:RBM2} of \(\boldsymbol X\), note that \(\REP_{\boldsymbol X}(\ak{\thetaN})\) has no immediate analytical expression similar to that of \(\REP_{(\boldsymbol X, \boldsymbol H)}(\ak{\thetaN})\). For clarity, recall also that S-instability \eqref{eq:Sun} in each model type refers to a respective divergence (i.e., \(\lim_{N\to \infty} \REP_{(\boldsymbol X, \boldsymbol H)}(\ak{\thetaN}) /(N+N_{\mathcal{H}})=\infty\), \(\lim_{N\to \infty} \REP(\ak{\thetaN}) /N=\infty\)) upon scaling by the corresponding number of variables in a distribution. In the following, let \(|\boldsymbol y|_1 = \sum_{i=1}^d |y_i|\) denote the L1 norm of a generic vector \(\boldsymbol y =(y_1,\ldots,y_d)\), \(d \geq 1\).
\BeginKnitrBlock{proposition}
\protect\hypertarget{prp:prop2}{}{\label{prp:prop2} }Let \(P_{\ak{\thetaN}}\) denote a RBM-based data model \eqref{eq:RBM2} for \(N\geq 1\) visible variables \(\boldsymbol X \equiv \boldsymbol X_N\) derived from \(\tilde{P}_{\ak{\thetaN}}\) as the joint RBM distribution \eqref{eq:RBM1} of \((\boldsymbol X, \boldsymbol H)\) involving some number \(N_{\mathcal{H}} \equiv N_{N,\mathcal{H}}\geq 0\) of hidden variables \(\boldsymbol H \equiv \boldsymbol H_N\) and parameters \(\ak{\thetaN} \equiv (\ak{\thetaN}^{\mathcal{V}},\ak{\thetaN}^{\mathcal{H}}, \ak{\thetaN}^{\mathcal{VH}}) \in\mathbb{R}^{N}\times \mathbb{R}^{N_{\mathcal{H}}} \times \mathbb{R}^{N*N_{\mathcal{H}}}\). Then,
\begin{enumerate}[(i)]
\item the instability measure $\REP(\ak{\thetaN})$ for the marginal model $P_{\ak{\thetaN}}$ of  $\boldsymbol X$ satisfies \ak{
    $$
    \left| \REP(\ak{\thetaN})  - \elt\right| \leq    N_{\mathcal{H}}  \log 2 \text{ for } \elt  \equiv   \max_{ \boldsymbol x} \max_{ \boldsymbol h  } f_{\ak{\thetaN}} (\boldsymbol x, \boldsymbol h)-\min_{ \boldsymbol x } \max_{ \boldsymbol h }f_{\ak{\thetaN}} (\boldsymbol x, \boldsymbol h) 
    $$}
    based on $f_{\ak{\thetaN}}$ from (\ref{eq:f}) with components $\boldsymbol x \in \{\pm 1\}^{N}, \boldsymbol h \in \{\pm 1\}^{N_{\mathcal{H}}}$.
\item The instability measure $\REP_{(\boldsymbol X, \boldsymbol H)}(\ak{\thetaN})\equiv \left(\max_{ \boldsymbol x} \max_{ \boldsymbol h  } f_{\ak{\thetaN}} (\boldsymbol x, \boldsymbol h)-\min_{ \boldsymbol x } \min_{ \boldsymbol h }f_{\ak{\thetaN}} (\boldsymbol x, \boldsymbol h)\right)$ for the joint model $\tilde{P}_{\ak{\thetaN}}$ of $(\boldsymbol X, \boldsymbol H)$ satisfies
    \begin{align*}
    2\Gam +  2|\ak{\thetaN}^{\mathcal{H}} |_{1} \; \geq \; \REP_{(\boldsymbol X, \boldsymbol H)}(\ak{\thetaN})    & \geq 
    2\max\big\{\Gam,   |\ak{\thetaN}^{\mathcal{H}} |_{1}\big\} \\
    &\geq  2\Gam\\
    &\geq  \elt\\
    &\geq   \max\big\{  \Gamc, \, \Gam  - 2|\ak{\thetaN}^{\mathcal{H}} |_{1}  \big\}
    \end{align*}
    for
    $$
    \Gam \equiv \max_{ \boldsymbol h} k_{\ak{\thetaN}} (\boldsymbol h)
    \geq  |\ak{\thetaN}^{\mathcal{V}} |_{1},\qquad  k_{\ak{\thetaN}} (\boldsymbol h)  \equiv \sum_{i=1}^{N }\left| \theta_{i,\ak{\thetaidx}}^{\mathcal{V}}   + \sum_{j=1}^{N_{\mathcal{H}}} h_j \theta_{ij,\ak{\thetaidx}}^{\mathcal{VH}} \right|,
    $$ 
    and $\Gamc \equiv \min_{ \boldsymbol h} k_{\ak{\thetaN}} (\boldsymbol h)$ based on a function  $k_{\ak{\thetaN}} (\boldsymbol h)$ of hidden variable outcomes $\boldsymbol h = (h_1,\ldots,h_{N_{\mathcal{H}}})$ and visible-related parameters $\ak{\thetaN}^{\mathcal{V}}$ and $\ak{\thetaN}^{\mathcal{VH}}$.
\item Assuming $\sup_{N \geq 1} N_{\mathcal{H}}/N<\infty$ additionally, then the following properties 1.-7. hold:
\begin{enumerate}[1.] \itemsep 0cm
    \item an S-unstable visible model $P_{\ak{\thetaN}}$ is equivalent to the condition $\lim_{N\to \infty}  \elt/N  =\infty$; further, $P_{\ak{\thetaN}}$ is stable when $\elt/N$, $N \geq 1$, is bounded.
    \item an S-unstable joint model $P_{\ak{\thetaN}}$ is equivalent to the condition $\lim_{N\to \infty} \max\{|\ak{\thetaN}^{\mathcal{H}}|_1, \Gam\}/N =\infty$; further,  $\tilde{P}_{\ak{\thetaN}}$ is stable when $[|\ak{\thetaN}^{\mathcal{H}}|_1+ \Gam]/N$, $N \geq 1$, is bounded.
    \item if the visible model $P_{\ak{\thetaN}}$   is S-unstable, then  the joint model $\tilde{P}_{\ak{\thetaN}}$ is also S-unstable.
    \item when $\lim_{N\to \infty}  (|\ak{\thetaN}^{\mathcal{V}}|_1-2|\ak{\thetaN}^{\mathcal{H}}|_1)/N =\infty$, both $P_{\ak{\thetaN}}$ and $\tilde{P}_{\ak{\thetaN}}$ are necessarily S-unstable.
    \item when $\lim_{N\to \infty}|\ak{\thetaN}^{\mathcal{H}}|_1/N =\infty$,  the joint model $\tilde{P}_{\ak{\thetaN}}$ is necessarily S-unstable.
    \item when $\sup_{N \geq 1} |\ak{\thetaN}^{\mathcal{H}}|_1 /N<\infty$, the visible model $P_{\ak{\thetaN}}$ being S-stable or S-unstable is equivalent to the joint model  $\tilde{P}_{\ak{\thetaN}}$ being stable or unstable.
    \item  an S-stable visible model $P_{\ak{\thetaN}}$ results if
        $$
        |\ak{\thetaN}^{\mathcal{V}}|_1+ |\ak{\thetaN}^{\mathcal{VH}} |_1  \leq CN,\quad N \geq 1,
        $$
        for some $C>0$, while an S-stable joint model $\tilde{P}_{\ak{\thetaN}}$  results if
        $$
        |\ak{\thetaN}|_1 \equiv  |\ak{\thetaN}^{\mathcal{V}}|_1+|\ak{\thetaN}^{\mathcal{H}}|_1 +|\ak{\thetaN}^{\mathcal{VH}} |_1  \leq  C N ,\quad N \geq 1.
        $$
\end{enumerate}
\end{enumerate}
\EndKnitrBlock{proposition}

\BeginKnitrBlock{remark}
{}The condition \(\sup_{N \geq 1} N_{\mathcal{H}}/N<\infty\) in Proposition \ref{prp:prop2}(iii) is often mild in practice (i.e., the number \(N_{\mathcal{H}}\) of hidden variables is typically not excessively larger than the number \(N\) of visible observations). This allows instability results for both marginal and joint RBM models to be more readily stated together, as the numbers \(N\) and \(N+N_{\mathcal{H}}\) of variables in these models become asymptotically equivalent.
\EndKnitrBlock{remark}

\ak{With regard to instability and the effects of different parameter types, the relationships between RBM models in Proposition \ref{prp:prop2}(iii)} follow from the bounds on model instability measures in Proposition \ref{prp:prop2}(i)-(ii). Generally speaking, all instability in the marginal RBM model for the data \(\boldsymbol X\) can be attributed to an excessively large model quantity \(\elt\), which predominantly follows when main \(\ak{\thetaN}^{\mathcal{V}}\) and interaction \(\ak{\thetaN}^{\mathcal{VH}}\) parameters related to visible variables are too large in magnitude (e.g., upon accumulation in terms such as \(|\ak{\thetaN}^{\mathcal{V}}|_1\), \(\Gam\) or \(\Gamc\)). For example, for any bounded sequence \(|\ak{\thetaN}^{\mathcal{H}}|/N\) of hidden parameters, \ak{main visible parameters $\ak{\thetaN}^{\mathcal{V}}$ that are too extreme ($|\ak{\thetaN}^{\mathcal{V}}|_1/N\to \infty$) will guarantee} instability in the visible model (\(\elt/N\to \infty\)). In fact, the instability measure \(\elt \equiv \max_{ \boldsymbol x} \max_{ \boldsymbol h } f_{\ak{\thetaN}} (\boldsymbol x, \boldsymbol h)-\min_{ \boldsymbol x } \max_{ \boldsymbol h }f_{\ak{\thetaN}} (\boldsymbol x, \boldsymbol h)\) for marginal/visible model represents a clearly smaller portion of the instability measure \(\REP_{(\boldsymbol X, \boldsymbol H)}(\ak{\thetaN})\equiv \max_{ \boldsymbol x} \max_{ \boldsymbol h } f_{\ak{\thetaN}} (\boldsymbol x, \boldsymbol h)-\min_{ \boldsymbol x } \min_{ \boldsymbol h } f_{\ak{\thetaN}} (\boldsymbol x, \boldsymbol h)\) in the joint RBM model\ak{. This implies} that an unstable marginal model (i.e., due to \(\ak{\thetaN}^{\mathcal{V}}\), \(\ak{\thetaN}^{\mathcal{VH}}\)) must always translate to an unstable joint model and that further potential causes of instability exist for the joint model, often due to the size \(|\ak{\thetaN}^{\mathcal{H}}|_1\).

\par

\ak{While} the joint RBM model for \((\boldsymbol X,\boldsymbol H)\) must always be unstable due to a diverging combination of visible and/or interaction parameters (\(|\ak{\thetaN}^{\mathcal{V}}|_1/N\to \infty\) or \(\Gam/N\to \infty\)) (Proposition \ref{prp:prop2}(iii.2)), instability for the joint model can also result when the main hidden parameters \(\ak{\thetaN}^{\mathcal{H}}\) become too large relative to sample size (\(|\ak{\thetaN}^{\mathcal{H}}|_1/N\to \infty\) in Proposition \ref{prp:prop2}(iii.5)). However, under Proposition \ref{prp:prop2}, the main hidden parameters \(\ak{\thetaN}^{\mathcal{H}}\) do not necessarily entail a source of instability for the marginal visible model\ak{. 
When} the hidden parameters are bounded relative to the sample size (\(\sup_{N\geq 1} |\ak{\thetaN}^{\mathcal{H}}|_1/N<\infty\)), then all instability in both the joint and marginal RBM models can be directly linked to excessively large visible \(\ak{\thetaN}^{\mathcal{V}}\) and/or interaction parameters \(\ak{\thetaN}^{\mathcal{VH}}\) so that features of stability/instability must be the same across both models (Proposition \ref{prp:prop2}(iii.6)). Hence, to prevent instability in the joint model, the combined magnitudes of all parameters \(\ak{\thetaN}\) must be controlled (cf.~Proposition \ref{prp:prop2}(iii.7)), while a stable visible data model technically results \ak{from} constraining only the sizes of visible-related parameters \(\ak{\thetaN}^{\mathcal{V}}\), \(\ak{\thetaN}^{\mathcal{VH}}\). Nevertheless, because the joint model \ak{often is} employed in practice for purposes of simulation and simulation-based inference, it is still reasonable to consider parameter choices for ensuring a stable joint model (and, consequently, a stable visible model as well). Further evidence of this is seen in the following numerical example.
\begin{figure}
\includegraphics{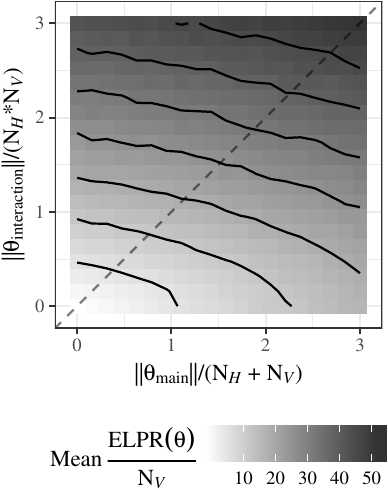} \includegraphics{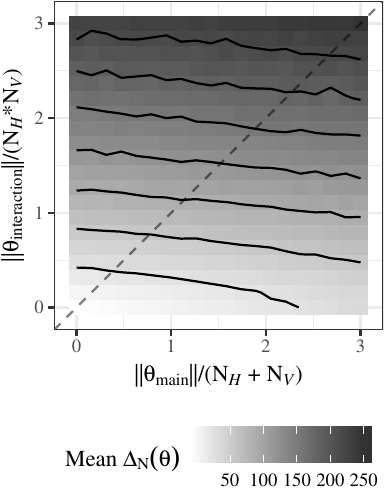} \caption{The sample mean value of $\text{ELPR}(\boldsymbol \theta)/N_{\mathcal{V}}$ (left) and $\DN(\boldsymbol \theta)$ at each grid point for each combination of magnitude of $\boldsymbol \theta$. As the magnitude of $\boldsymbol \theta$ grows, so does the value of these metrics, indicating \ak{instability} in the model.}\label{fig:rbm-plots}
\end{figure}
\par

In our numerical experiment, we allow the two types of terms (main effects terms corresponding to visible and hidden parameters \(\boldsymbol \theta_{main} = (\ak{\thetaN}^{\mathcal{V}}, \ak{\thetaN}^{\mathcal{H}})\) and interaction parameters \(\ak{\thetaN}^{\mathcal{VH}}\)) to have varying average magnitudes, \(||\boldsymbol \theta_{main} || /(N_{\mathcal{H}}+N_{\mathcal{V}})\) and \(||\boldsymbol \theta_{interaction} || /(N_{\mathcal{H}}*N_{\mathcal{V}})\) for a RBM with \(N_\mathcal{V} = 9\) visibles and \(N_\mathcal{H} = 5\) hiddens. These average magnitudes vary on a grid between \(0.001\) and \(3\) with \(20\) breaks, yielding \(400\) grid points. At each point in the grid, \(100\) vectors (\(\boldsymbol \theta_{main}\)) are sampled uniformly on a sphere with radius corresponding to the first coordinate in the grid and \(100\) vectors (\(\boldsymbol \theta_{interaction}\)) are sampled uniformly on a sphere with radius corresponding to the second coordinate in the grid via sums of squared and scaled iid Normal\((0, 1)\) variables. These vectors are then paired to create \(100\) values of \(\ak{\thetaN}\) with magnitudes at each point in the grid. The values \(\REP(\ak{\thetaN})/N_{\mathcal{V}}\) and \(\DN(\ak{\thetaN})\) are then calculated for each \(\ak{\thetaN}\) and then summarized for each point in the grid using the sample mean. The results of this numerical study are shown in Figure \ref{fig:rbm-plots}. From these two plots, it is clear that for larger magnitudes of the parameter vectors, there is evidence of S-instability in that the log-ratio of extremal probabilities scaled by \(N_{\mathcal{V}}\) and the the biggest log-probability ratio for a one-component change in data outcomes are both increasing away from \(\ak{\thetaN} = \boldsymbol 0\), further supporting \ref{prp:prop2}(iii.2 and iii.5).

\par

In more complicated graphical models involving further or deeper hidden layers, the same issues and causes of instability similarly exist, but are compounded by a greater number of model parameters \ak{and depend on the configuration of the particular models in consideration. Consider the following two models -- a RBM with $N$ visible nodes and one hidden layer with $N_H$ nodes and a deep RBM with $N$ visible nodes and $2$ layers, with $N_{H,1}$ and $N_{H,2}$ hidden nodes respectively such that $N_{H,1} + N_{H,2} = N_H$. The length of the parameter vector $\ak{\thetaN}$ in the single layer RBM is then equal to $N + N_H + N*N_H$, whereas the length of the parameter vector in the two layer RBM is $N + N_{H,1} + N_{H,2} + N*N_{H,1} + N_{H,1}*N_{H,2} = N + N_H + N*N_H + N_{H,2}(N_{H,1} - N)$. If $N_{H,1} > N$, then the deep model will have more parameters, if $N_{H,1} < N$, then the deep model will have less parameters, and if $N_{H,1} = N$, then the deep model will the same number of parameters as the single layer RBM.} S-unstable joint models will similarly follow if the combined sizes of all parameters are too great relative to the total number of variables, while instability in the data model for visible variables will depend only on the main or interaction parameters directly related to visibles \ak{(i.e. all the parameters for the single layer RBM and only the parameters for the visible layer, layer $H^{(1)}$, and the interactions between them for the two layer RBM)} and how their accumulated magnitude compares to the observation sample size \(N\). \ak{This indicates that for the name number of hidden variables $N_H$, a deeper, rather than a wider, structure can be beneficial when considering instability in the data model for visible variables. However, as the size of a model considered grows, this point becomes moot.}

\hypertarget{implications}{%
\section{Statistical Consequences of Instability}\label{implications}}

Due to \ak{sensitivity in the probability structure of an unstable FOES model, S-instability} may often translate to numerical complications, and in fact obstructions, in both simulation and statistical inference based on likelihoods. We describe these aspects in Sections \ref{ml}-\ref{bayes} with regard to data simulation, maximum likelihood estimation and Bayes inference.

\hypertarget{ml}{%
\subsection{Implications for Maximum Likelihood Inference}\label{ml}}

Volatility in the probability structure of an unstable model can also hamper efforts to maximize likelihood functions in statistical inference. When a FOES model is unstable along a parameter sequence \(\ak{\thetaN}\), the same model can further be unstable along parameters \(-\ak{\thetaN}\) in an opposite direction from the origin (model condition PSR and Corollary \ref{cor:sign}). This can translate into potential sensitivity of the likelihood function around zero, and lead to numerical complications in maximizing the objective function. We next provide a discussion of this issue in a way that builds upon and extends related findings by \citet{schweinberger2011instability}, who largely focused on the case of one-parameter exponential models.

With many probability models, the modes and anti-modes in the probability structure under one parameter \(\ak{\thetaN}\) are reversed in role when the parameter sign changes \(-\ak{\thetaN}\). Because unstable models tend to degeneracy, the opposite signed parameters further push unstable models to assign nearly all probability to extremely opposite data configurations, given by modes/anti-modes. This is made concrete in Theorem \ref{thm:something}, relating the degeneracy from unstable models to the expected behavior of log-likelihood functions.

\BeginKnitrBlock{theorem}
\protect\hypertarget{thm:something}{}{\label{thm:something} }Let \(P_{\ak{\thetaN}}\), \(N \geq 1\), denote an S-unstable FOES model sequence, which additionally satisfies model condition PSR. Let \(\boldsymbol{x}_{\max, \ak{\thetaN}},\boldsymbol{x}_{\min, \ak{\thetaN}}\in\mathcal{X}^N\) denote, respectively, a mode and anti-mode of the model \(P_{\ak{\thetaN}}(\boldsymbol x)\), \(\boldsymbol x\in\mathcal{X}^N\), for \(N\) observations \(\boldsymbol X = (X_1,\ldots,X_N)\), whereby \(P_{\ak{\thetaN} }(\boldsymbol{x}_{\max, \ak{\thetaN}}) = \max_{\boldsymbol y\in\mathcal{X}^N} P_{\ak{\thetaN} }(\boldsymbol y)\) and \(P_{\ak{\thetaN} }(\boldsymbol{x}_{\min, \ak{\thetaN}}) = \min_{\boldsymbol y\in\mathcal{X}^N} P_{\ak{\thetaN} }(\boldsymbol y)\).

Then, letting \(\stackrel{p, \mathrm{E}}{\longrightarrow}\) denote convergence in probability and expectation, as \(N\to \infty\),
\[
\frac{1}{\REP(\ak{\thetaN})} \log\left[ \frac{P_{\ak{\thetaN} }(\boldsymbol X) }{ \displaystyle{\min_{\boldsymbol y\in\mathcal{X}^N}} P_{\ak{\thetaN} }(\boldsymbol y)}\right] = \frac{ \log P_{\ak{\thetaN} } ( \boldsymbol X)  - \log P_{\ak{\thetaN} }(\boldsymbol{x}_{\min, \ak{\thetaN}}) }
{ \log P_{\ak{\thetaN} } (\boldsymbol{x}_{\max,\ak{\thetaN}})  - \log P_{\ak{\thetaN} }(\boldsymbol{x}_{\min, \ak{\thetaN}}) }
\stackrel{p, \mathrm{E}}{\longrightarrow}
\; 1
\]
under \(\ak{\thetaN}\) while
\[
\frac{1}{\REP(-\ak{\thetaN})} \log\left[ \frac{P_{-\ak{\thetaN} }(\boldsymbol X) }{  \displaystyle{\min_{\boldsymbol y\in\mathcal{X}^N}} P_{-\ak{\thetaN} }(\boldsymbol y)} \right]
=    \frac{ \log P_{-\ak{\thetaN} } ( \boldsymbol X) - \log P_{\ak{\thetaN} }(\boldsymbol{x}_{\max, \ak{\thetaN}})    }
{ \log P_{-\ak{\thetaN} } (\boldsymbol{x}_{\min, \ak{\thetaN}})  - \log P_{-\ak{\thetaN} }(\boldsymbol{x}_{\max, \ak{\thetaN}}) }
\stackrel{p,\mathrm{E}}{\longrightarrow}
\; 1,
\]
under \(-\ak{\thetaN}\), where
\[
\REP( \ak{\thetaN}) \equiv \log \frac{P_{\ak{\thetaN} } (\boldsymbol{x}_{\max, \ak{\thetaN}})}{P_{\ak{\thetaN} }(\boldsymbol{x}_{\min, \ak{\thetaN}})} = \log  \frac{P_{-\ak{\thetaN} } (\boldsymbol{x}_{\min, \ak{\thetaN}}) }{ P_{-\ak{\thetaN} }(\boldsymbol{x}_{\max, \ak{\thetaN}})}=   \REP( -\ak{\thetaN}), \quad N \geq 1.
\]
\EndKnitrBlock{theorem}

Theorem \ref{thm:something} \ak{implies} log-likelihood functions based on unstable models are both inversely related and degenerate at opposited signed parameters \(\ak{\thetaN}\) or \(-\ak{\thetaN}\), so that likelihoods are highest at different extremes in data configuration (e.g., \(\boldsymbol{x}_{\max, \ak{\thetaN}}\) under \(\boldsymbol \theta\)-probabilities or \(\boldsymbol{x}_{\min, \ak{\thetaN}}\) under \(-\boldsymbol \theta\)-probabilities). If the observed outcome \(\boldsymbol x\) for data \(\boldsymbol X\) is not a mode/anti-mode, then probabilities for the outcome may be small under both parameters \(\ak{\thetaN}\) and \(-\ak{\thetaN}\), in which case associated optimization steps may then shift around zero and struggle to converge.

In many model formulations, the zero parameter \(\ak{\thetaN}=\boldsymbol 0\) is a ``safe'' position among parameters, representing a guaranteed stable model (having uniform probabilities among outcomes), which can also tether a broad parameter search attempted among unstable models. \citet{handcock2003assessing} describes similar results for degenerate exponential models, and Theorem \ref{thm:something} also supports an important finding of \citet[Corollary 1]{schweinberger2011instability} for one-parameter linear exponential models \eqref{eq:expo}. In the latter case, the likelihood score function at \(\ak{\thetaN}\) is the expected value \(\mu( \ak{\thetaN})\equiv \mathrm{E}_{\ak{\thetaN}} g(\boldsymbol X)\) of the sufficient statistic \(g(\cdot)\), and optimization involves solving \(\mu(\cdot)=g(\boldsymbol x)\) for an observed outcome \(\boldsymbol x\). For unstable models in this exponential class, \citet[Corollary 1]{schweinberger2011instability} shows
that
\[
\lim_{N\to \infty}\frac{\mu(\ak{\thetaN}) -L_n}{U_n-L_N}= \left\{ \begin{array}{lcl}
1 && \text{for } \ak{\thetaN}>0,\\
0 && \text{for } \ak{\thetaN}<0, \end{array}\right.
\]
where again \(U_N\) and \(L_N\) denote the maximum and minimum values of the statistic \(g(\boldsymbol x)\), \(\boldsymbol x\in\mathcal{X}^N\). As described by \citet{schweinberger2011instability}, the implication for maximum likelihood estimation is that, unless an observed outcome \(\boldsymbol x\) falls at an extreme \(U_N, L_N\) (i.e., modes/anti-modes), optimization steps in the parameter space can iterate in relatively small increments around zero and fail to converge. For unstable one-parameter exponential models, the maximum likelihood results of \citet{schweinberger2011instability} turn out to be a special case of Theorem \ref{thm:something} and the LREP expansion \eqref{eq:UL} in this setting; namely, for an unstable model with \(\ak{\thetaN}>0\),
\[
\frac{1}{\REP(\ak{\thetaN})} \log\left[ \frac{P_{\ak{\thetaN} }(\boldsymbol X) }{\displaystyle{\min_{\boldsymbol y\in\mathcal{X}^N}} P_{ \ak{\thetaN} }(\boldsymbol y)}\right] = \frac{\boldsymbol g(\boldsymbol X) - L_N}{U_N-L_N}\;\;\stackrel{p, \mathrm{E}}{\longrightarrow}\;\; 1
\]
holds as \(N\to \infty\) by Theorem \ref{thm:something}, while under \(-\ak{\thetaN}<0\)
\[
\frac{1}{\REP(-\ak{\thetaN})} \log\left[ \frac{P_{-\ak{\thetaN} }(\boldsymbol X) }{\displaystyle{\min_{\boldsymbol y\in\mathcal{X}^N}} P_{-\ak{\thetaN} }(\boldsymbol y)}\right] = \frac{U_N -\boldsymbol g(\boldsymbol X)}{U_N-L_N} = 1 - \frac{\boldsymbol g(\boldsymbol X) - L_N}{U_N-L_N} \;\;\stackrel{p, \mathrm{E}}{\longrightarrow}\;\; 1.
\]
Again, when all probability in unstable models may be pushed to opposite extremes in the sample space, due to a combination of degeneracy and parameter sign, numerical complications in likelihood maximization may occur.

\hypertarget{bayes}{%
\subsection{Implications for Bayes Inference}\label{bayes}}

The potential numerical difficulties with maximum likelihood with unstable models, as described in the previous section, can naturally carry over to Bayes inference. Considering that the degeneracy issues related to unstable models can cause likelihoods can be flat (e.g., near zero) for many parameters under a given data outcome and that sign changes in parameters can shift tremendous probability to extreme and opposite outcomes in the sample space (e.g., Corollary \ref{cor:sign}, Theorem \ref{thm:something}), then numerical complications may arise with Bayes inference in sampling a posterior parameter space based on MCMC.
\ak{Instability may hinder effective chain mixing, related to the findings of \citet{handcock2003assessing} and \citet{schweinberger2011instability} with exponential graph models. For example, a Markov chain may become entrapped within a mode, and fail to adequately mimic the overall occupation frequencies required for a reasonable MCMC distributional approximaton. If modes of the unstable model are not unique, then important outcomes may be missed without multiple chains or impractically enormous numbers of MCMC samples. This mixing problem is due to the unstable stationary distribution (unbounded ratios of probabilities under the joint model), rather than in any particulars of the MCMC algorithm. Hence,} in the Bayes setting for sampling a posterior distribution for \(\ak{\thetaN}\), a chain may \ak{be unable} to effectively explore the parameter space due partly to extreme and potentially unbounded probability ratios from parameter sign changes.

For example, if \(\pi(\cdot)\) denotes a prior density for \(\ak{\thetaN}\) and \(q(\cdot | \cdot)\) denotes a proposal distribution for use in a Metropolis-Hastings (MH) sampler, then MH acceptance probability becomes
\[
\alpha\left(\boldsymbol \theta^{(1)} \mid \boldsymbol \theta^{(2)}\right)= \min\left\{1,
\frac{q(\ak{\thetaN}^{(2)} \mid \ak{\thetaN}^{(1)})}{q(\ak{\thetaN}^{(1)} \mid  \ak{\thetaN}^{(2)})}
\frac{P_{\ak{\thetaN}^{(1)}} ( \boldsymbol x ) \pi(\ak{\thetaN}^{(1)}) }{P_{\ak{\thetaN}^{(2)}} (  \boldsymbol x) \pi(\ak{\thetaN}^{(2)}) } \right\},
\]
which indicates how parameter sensitivity in the likelihood \(P_{\ak{\thetaN}} (\boldsymbol x)\) may complicate sampling of the posterior \(P_{ \ak{\thetaN}} ( \boldsymbol x) \pi(\ak{\thetaN})\) (i.e., moving from \(\boldsymbol \theta^{(1)}\) to \(\boldsymbol \theta^{(2)}\) in the parameter space). Furthermore, the potential for model instability and the size of the parameter space can also become greater with the introduction of latent variables to existing data variables, as involved in some model formulations described in Sections \ref{rbm}-\ref{deep-learning}. As latent variables are often sampled with parameters in a Bayes MCMC approach, this aspect may further compound numerical problems in chain mixing.

\hypertarget{stable-models}{%
\subsection{Fitting Stable Models}\label{stable-models}}

\ak{The instability results presented in Section \ref{instability-results} also suggest the potential use of regularization and penalization as a solution to avoid instability when fitting FOES models. One example of a rigorous approach to penalization can be found in \citet{kaplan2019properties}, which employs Bayesian fitting for a single layer RBM using multivariate Gaussian priors with constrained covariance matrices to control the sizes of the parameter values. For example, considering a latent variable model described in Section \ref{fixed-dim-exp}, Proposition \ref{prp:prop2}(iii) suggests how parameter contributions from sources of ``main hidden," ``main visible," and ``interaction" parameters may each be constrained (relative to the visible variable sample size $N$) to ensure model stability. In a Bayes inference framework, priors can be developed so that parameters are then appropriately bounded (in a stochastic sense) to lie in a parameter subset compatible with stability. For example, for a latent model with numbers $N_{\mathcal{H}}$, $N_{\mathcal{V}}$ and $N_{\mathcal{VH}}$ for hidden, visible and interaction parameters, respectively, Proposition \ref{prp:prop2}(iii.7) suggests a potential prior specification where parameters have root second moments bounded by the parameter type: $C N/N_{\mathcal{V}}$, $C N/N_{\mathcal{H}}$, and $C N/N_{\mathcal{VH}}$ for visible, hidden and interactions parameters with a constant $C>0$ (i.e., so that the condition in Proposition \ref{prp:prop2}(iii.7) holds in a stochastic sense). The same general strategy can apply with other models with results from Section \ref{illustrations}, and the adequacy of the subsequent fitted model might then be assessed.

As opposed to approaches for fitting (e.g., penalization as above),  the findings here on model stability/instability also support alternative strategies to model formulation and re-pameterization using ``centered conditional" distributions \citep[cf.][]{kaiser2007statistical}; see \citet{kaiser2009exploring} and \citet{kaiser2012centered} for applications with spatial data and \citet{casleton2017local} for network data. In a joint model $P_{\thetaN}$,  these authors consider associated (full) conditional distributions $P_{\thetaN}(X_{i} = x_i | X_j=x_j, j\neq i)$ and re-parameterize by centering/averaging sums in the expression of  $P_{\thetaN}(X_{i} = x_i | X_j=x_j, j\neq i)$ in order to separate model effects between ``mean" and ``dependence" parameters.  The intent is to improve model interpretation and help detect degeneracy (characterized by the ``mean parameters" in the model failing to match averages observed in data generations from the model and attributable to large ``dependence parameters"). The works above describe advantages of such centered parameterizations based on empirical studies, but findings here also support these parameterizations as useful for avoiding model instability, where the corresponding ``centered" models can be checked for $\REP(\thetaN)/N$ ratios that are bounded over wider possible parameters than ``uncentered" counterpart models.  For illustration, the appendix provides details on a graph model example.}

\hypertarget{conclusions}{%
\section{Concluding Remarks}\label{conclusions}}

For a large class of models that covers a broad range of applications (including ``deep learning''), we have developed a formal definition of instability in model probability structure and elucidated multiple consequences of instability. We have shown for FOES models that instability manifests through \ak{a large amount of probability being placed on a subset of the theoretically possible realizations that may be narrower than intended.}
Such instability is often due to a complex interaction between the model statistics used (i.e., how numerous and large these may become) and the number and magnitudes of parameters in the model formulation. For many FOES models, the possibility exists, at least in principle, to \ak{constrain} parameters in a way \ak{that} balances their potential contributions against those of model statistics in order to prevent probability instabilities. \ak{While the results on instability here only address model formulations with finite sample spaces, the same underlying principles might be extended to models for continuous data as well by discretizing the model (i.e., binning outcomes into intervals), as data typically have some given level of discretization. Such a generalization could help toward understanding potential model instability in a wider class of data applications. As it sits currently, t}he FOES model class is quite broad and \ak{findings here can help in identifying undesirable probability features in such models.}

\clearpage

\hypertarget{appendix-appendix}{%
\appendix}

\hypertarget{appendix-instab}{%
\section{Proofs of instability results}\label{appendix-instab}}

\textbf{Proof of Theorem \ref{thm:instab-elpr}.} For part (i), we prove the contrapositive, supposing that \(\DN(\ak{\thetaN}) \le C\) holds for some \(C > 0\) and show \(\REP(\ak{\thetaN}) \leq NC\). Let \(\boldsymbol x_{min} \equiv \argmin\limits_{\boldsymbol x \in \mathcal{X}^N}P_{\ak{\thetaN}}(\boldsymbol x)\)
and \(\boldsymbol x_{max} \equiv \argmax\limits_{\boldsymbol x \in \mathcal{X}^N}P_{\ak{\thetaN}}(\boldsymbol x)\). Note there exists a sequence \(\boldsymbol x_{min} \equiv \boldsymbol x_0, \boldsymbol x_1, \dots, \boldsymbol x_k \equiv \boldsymbol x_{max}\) in \(\mathcal{X}^N\) of component-wise switches to move from \(\boldsymbol x_{min}\) to \(\boldsymbol x_{max}\) in the sample space (i.e.~\(\boldsymbol x_i, \boldsymbol x_{i + 1} \in \mathcal{X}^N\) differ in exactly \(1\) component, \(i = 0, \dots, k\)) for some integer \(k \in \{0, 1, \dots, N\}\). Under the FOES model, recall
\(P_{\ak{\thetaN}}(\boldsymbol x) > 0\) holds so that \(\log P_{\ak{\thetaN}}(\boldsymbol x)\) is well-defined for each outcome \(\boldsymbol x \in \mathcal{X}^N\). Then, if \(k > 0\), it follows that
\begin{align*}
\REP(\ak{\thetaN}) = \log\left[\frac{P_{\ak{\thetaN}}(\boldsymbol x_{max})}{P_{ \ak{\thetaN}}(\boldsymbol x_{min})}\right] &= \left|\sum\limits_{i = 1}^k\log\left(\frac{P_{\ak{\thetaN}}(\boldsymbol x_i)}{P_{\ak{\thetaN}}(\boldsymbol x_{i-1})}\right)\right| \\
&\le \sum\limits_{i = 1}^k\left|\log\left(\frac{P_{\ak{\thetaN}}(\boldsymbol x_i)}{P_{\boldsymbol \theta}(\boldsymbol x_{i-1})}\right)\right| \le k \Delta_N(\ak{\thetaN}) \le NC,
\end{align*}
using \(k \le N\) and \(\Delta(\ak{\thetaN}) \le C\). If \(k = 0\), then \(\boldsymbol x_{max} = \boldsymbol x_{min}\) and the same bound above holds. This establishes part (i). To show part (ii), note the definition of S-instability (i.e., \(\lim_{N\to \infty}\REP(\ak{\thetaN})/N= \infty\)) combined with part (i) implies that \(\lim_{N\to \infty}\DN(\ak{\thetaN})=\infty\). \hfill \(\Box\)

\textbf{Proof of Theorem \ref{thm:degenFOES}.} As \(|\mathcal{X}|<\infty\) holds in the FOES model, we may suppose
\(|\mathcal{X}|>1\); otherwise, \(\mathcal{X}^N\) has one outcome and the model is trivially degenerate for all \(N \geq 1\). Fix \(0 < \epsilon < 1\) and write \(\boldsymbol x_{min} \equiv \argmin\limits_{\boldsymbol x \in \mathcal{X}^N}P_{ \ak{\thetaN}}(\boldsymbol x)\) and \(\boldsymbol x_{max} \equiv \argmax\limits_{\boldsymbol x \in \mathcal{X}^N}P_{\ak{\thetaN}}(\boldsymbol x)\). Then, \(\boldsymbol x_{max} \in M_{\epsilon, \ak{\thetaN}}\), so \(P_{\ak{\thetaN}}(M_{\epsilon, \ak{\thetaN}}) \ge P_{ \ak{\thetaN}}(\boldsymbol x_{max}) > 0\). If \(\boldsymbol x \in \mathcal{X}^N \setminus M_{\epsilon, \ak{\thetaN}}\), then by definition
\(P_{\ak{\thetaN}}(\boldsymbol x) \le [P_{\ak{\thetaN}}(\boldsymbol x_{max})]^{1-\epsilon}[P_{\ak{\thetaN}}(\boldsymbol x_{min})]^{\epsilon}\) holds so that
\[
1-P_{\ak{\thetaN}}(M_{\epsilon, \ak{\thetaN}})
 = \sum\limits_{\boldsymbol x \in \mathcal{X}^N \setminus M_{\epsilon, \ak{\thetaN}}}P_{\ak{\thetaN}}(\boldsymbol x)
  \le (|\mathcal{X}|^N)[P_{\ak{\thetaN}}(\boldsymbol x_{max})]^{1-\epsilon}[P_{ \ak{\thetaN}}(\boldsymbol x_{min})]^{\epsilon}.
\]
From the lower bound on \(P_{\ak{\thetaN}}(M_{\epsilon, \ak{\thetaN}})\) and the
upper bound on \(1-P_{\ak{\thetaN}}(M_{\epsilon, \ak{\thetaN}})\), it follows that
\begin{align*}
\frac{1}{N}\log\left[\frac{P_{\ak{\thetaN}}(M_{\epsilon, \ak{\thetaN}})}{1-P_{ \ak{\thetaN}}(M_{\epsilon, \ak{\thetaN}})}\right] & \ge \frac{1}{N} \log\left[\frac{P_{\ak{\thetaN}}(\boldsymbol x_{max})}{(|\mathcal{X}|^N)[P_{ \ak{\thetaN}}(\boldsymbol x_{max})]^{1-\epsilon}[P_{\ak{\thetaN}}(\boldsymbol x_{min})]^{\epsilon}}\right] \\
&= \frac{\epsilon}{N} \log\left[\frac{P_{\ak{\thetaN}}(\boldsymbol x_{max})}{P_{ \ak{\thetaN}}(\boldsymbol x_{min})}\right] - \log |\mathcal{X}| \rightarrow \infty
\end{align*}
as \(N \rightarrow \infty\) by the definition of an S-unstable FOES model \eqref{eq:Sun}. Consequently, \(P_{\ak{\thetaN}}(M_{\epsilon, \ak{\thetaN}}) \rightarrow 1\) as \(N \rightarrow \infty\) as claimed. \hfill \(\Box\)

\textbf{Proof of Corollary \ref{cor:sign}.} The model condition PSR implies that
\begin{equation}
\label{eq:R}
\frac{\max_{\boldsymbol y \in \mathcal{X}^N}P_{\ak{\thetaN}}(\boldsymbol y)  }{\min_{\boldsymbol y \in \mathcal{X}^N}P_{\ak{\thetaN}}(\boldsymbol y) } = \frac{\max_{\boldsymbol y \in \mathcal{X}^N}P_{-\ak{\thetaN}}(\boldsymbol y) }{\min_{\boldsymbol y \in \mathcal{X}^N}P_{-\ak{\thetaN}}(\boldsymbol y) }
\end{equation}
so that the log-ratio \(\REP(\ak{\thetaN})=\REP(-\ak{\thetaN})\) is the same for both \(\ak{\thetaN}\) and \(-\ak{\thetaN}\) in \eqref{eq:elpr}. Now part (i) of Corollary \ref{cor:sign} follows from \(\REP(\ak{\thetaN})/N=\REP(-\ak{\thetaN})/N\to \infty\) as \(N\to \infty\) in \eqref{eq:Sun}. To show part (ii), fix \(0 < \epsilon < 1\) and consider a \(\epsilon\)-mode set \(\mathcal{M}_{\epsilon, \ak{\thetaN}}\) under \(\ak{\thetaN}\) from \eqref{eq:mode}. If \(\boldsymbol x \in \mathcal{M}_{\epsilon, \ak{\thetaN}}^c \equiv \mathcal{X}^N \setminus \mathcal{M}_{\epsilon, \ak{\thetaN}}\), then, by definition,
\[
\frac{P_{\ak{\thetaN}}(\boldsymbol x)}{\min_{\boldsymbol y \in \mathcal{X}^N}P_{ \ak{\thetaN}}(\boldsymbol y)} \leq \left[\frac{\max_{\boldsymbol y \in \mathcal{X}^N}P_{\ak{\thetaN}}(\boldsymbol y)  }{\min_{\boldsymbol y \in \mathcal{X}^N}P_{\ak{\thetaN}}(\boldsymbol y) } \right]^{1-\epsilon}
\]
holds, which is equivalent to
\[
\frac{\max_{\boldsymbol y \in \mathcal{X}^N}P_{-\ak{\thetaN}}(\boldsymbol y)}{P_{-\ak{\thetaN}}(\boldsymbol x)} \leq
\left[\frac{\max_{\boldsymbol y \in \mathcal{X}^N}P_{-\ak{\thetaN}}(\boldsymbol y)  }{\min_{\boldsymbol y \in \mathcal{X}^N}P_{-\ak{\thetaN}}(\boldsymbol y) } \right]^{1-\epsilon}
\]
by model condition PSR and \eqref{eq:R}. The latter is in turn equivalent to
\begin{equation}
\label{eq:R2}
\log P_{-\ak{\thetaN}}(\boldsymbol x) \geq \epsilon \max\limits_{\boldsymbol y \in \mathcal{X}^N} \log  P_{-\ak{\thetaN}}(\boldsymbol y) + (1-\epsilon)\min\limits_{\boldsymbol y \in \mathcal{X}^N} \log P_{-\ak{\thetaN}}(\boldsymbol y),
\end{equation}
so that \(\boldsymbol x \in \mathcal{M}_{\epsilon, \ak{\thetaN}}^c\) if and only if \eqref{eq:R2} holds. Next consider the \((1-\epsilon)\)-mode set \(\mathcal{M}_{1-\epsilon, -\ak{\thetaN}}\) under \(-\ak{\thetaN}\) from \eqref{eq:mode}. If \(\boldsymbol x \in\mathcal{M}_{1-\epsilon, -\ak{\thetaN}}\), then by definition \(\boldsymbol x\) fulfills \eqref{eq:R2} and so \(\boldsymbol x \in \mathcal{M}_{\epsilon, \ak{\thetaN}}^c\), showing that \(\mathcal{M}_{1-\epsilon, -\ak{\thetaN}} \subset \mathcal{M}_{\epsilon, \ak{\thetaN}}^c\). By this and the fact that that Theorem \ref{thm:degenFOES} and Corollary \ref{cor:sign}(i) entail that \(P_{-\ak{\thetaN}}(\boldsymbol X_N \in \mathcal{M}_{1-\epsilon, -\ak{\thetaN}})\to 1\) as \(N\to \infty\) (i.e., \(P_{-\ak{\thetaN}}\) is S-unstable), we have
\[
1  = \lim_{N\to \infty} P_{-\ak{\thetaN}}(\boldsymbol X_N \in \mathcal{M}_{1-\epsilon, -\ak{\thetaN}}) \leq \lim_{N\to \infty} P_{-\ak{\thetaN}}(\boldsymbol X_N \in \mathcal{M}_{\epsilon, \ak{\thetaN}}^c ) \leq 1,
\]
proving Corollary \ref{cor:sign}(ii) \hfill \(\Box\)

\textbf{Proof of Proposition \ref{prp:prop1}.} For any two outcomes \(\boldsymbol x_1, \boldsymbol x_2\in\mathcal{X}^N\), the log-ratio of probabilities from the linear exponential model \eqref{eq:expo} with \(k\) parameters/statistics satisfies
\[
\left|\log \left[ \frac{P_{\ak{\thetaN}}(\boldsymbol x_1)}{P_{\ak{\thetaN}}(\boldsymbol x_2)}  \right] \right| =
\left|  \sum_{i=1}^k \theta_{i,\ak{\thetaidx}} [g_{i,k}(\boldsymbol x_1) - g_{i,k}(\boldsymbol x_2) ] \right|  \leq  \sum_{i=1}^k | \theta_{i,\ak{\thetaidx}}| (U_{i,N}-L_{i,N});
\]
consequently, \(\REP(\ak{\thetaN} ) \leq \sum_{i=1}^k | \theta_{i,\ak{\thetaidx}}| (U_{i,N}-L_{i,N})\) holds in \eqref{eq:elpr} and model stability in Proposition \ref{prp:prop1} follows from \eqref{eq:Sun}. \hfill \(\Box\)

\textbf{Proof of Proposition \ref{prp:prop2}.}
Writing \(\boldsymbol x=(x_1,\ldots,x_N)\) and \(\boldsymbol h = (h_1,\ldots,h_{N_{\mathcal{H}}})\) with all components \(x_i,h_j\in\{\pm 1\}\), probabilities in the joint RBM model \eqref{eq:RBM1} can be written as \(\tilde{P}_{\ak{\thetaN}} (\boldsymbol x, \boldsymbol h) = c(\ak{\thetaN})\exp[ f_{\ak{\thetaN}} (\boldsymbol x, \boldsymbol h)]\) in terms of the function \(f_{\ak{\thetaN}} (\boldsymbol x, \boldsymbol h)\) from \eqref{eq:f} and the normalizing constant \(c(\ak{\thetaN})= \exp [-\psi(\ak{\thetaN})]\) from \eqref{eq:RBM1}. Let \(\boldsymbol x_M, \boldsymbol x_m\in\{\pm 1\}^N\) be such that \(P_{\ak{\thetaN}} (\boldsymbol x_M) = \max_{\boldsymbol x}P_{\ak{\thetaN}} (\boldsymbol x)\) and \(P_{\ak{\thetaN}} (\boldsymbol x_m) = \min_{\boldsymbol x}P_{\ak{\thetaN}} (\boldsymbol x)\) under the marginal RBM model \(P_{\ak{\thetaN}} (\boldsymbol x) = c(\ak{\thetaN})\sum_{\boldsymbol h \in\{\pm 1\}^{\mathcal{N}_H}} \tilde{P}_{\ak{\thetaN}} (\boldsymbol x, \boldsymbol h)= c( \ak{\thetaN})\sum_{\boldsymbol h \in\{\pm 1\}^{\mathcal{N}_H}} \exp[ f_{\ak{\thetaN}} (\boldsymbol x, \boldsymbol h)]\) from \eqref{eq:RBM2}. Also, \(\boldsymbol x_0,x_1\in\{\pm 1\}^N\) be such that \(\max_{\boldsymbol h}f_{\ak{\thetaN}} (\boldsymbol x_0, \boldsymbol h)=\max_{\boldsymbol x}\max_{\boldsymbol h}f_{\ak{\thetaN}} (\boldsymbol x , \boldsymbol h)\) and \(\max_{\boldsymbol h}f_{\ak{\thetaN}} (\boldsymbol x_1, \boldsymbol h)=\min_{\boldsymbol x}\max_{\boldsymbol h}f_{\ak{\thetaN}} (\boldsymbol x , \boldsymbol h)\). Then, Proposition \ref{prp:prop2}(i) follows from \(\REP_{(\boldsymbol X)}(\ak{\thetaN}) = \log[P_{\ak{\thetaN}} (\boldsymbol x_M) /P_{\ak{\thetaN}} (\boldsymbol x_m) ]\) and the lower/upper bounds on \(P_{\ak{\thetaN}} (\boldsymbol x_M)\) and \(P_{\ak{\thetaN}} (\boldsymbol x_m)\) as
\[
c(\ak{\thetaN}) \exp[\max_{\boldsymbol h}f_{\ak{\thetaN}} (\boldsymbol x_0 , \boldsymbol h)]
\leq P_{\ak{\thetaN}} (\boldsymbol x_0) \leq  P_{\ak{\thetaN}} (\boldsymbol x_M) \leq  2^{N_{\mathcal{H}}} c(\ak{\thetaN}) \exp[\max_{\boldsymbol x}\max_{\boldsymbol h}f_{\ak{\thetaN}} (\boldsymbol x, \boldsymbol h)]
\]
and
\begin{align*}
c(\ak{\thetaN}) \exp[\min_{\boldsymbol x}\max_{\boldsymbol h}f_{\ak{\thetaN}} (\boldsymbol x , \boldsymbol h)] &\leq c(\ak{\thetaN}) \exp[ \max_{\boldsymbol h}f_{\ak{\thetaN}} (\boldsymbol x_m , \boldsymbol h)] \\
&\leq P_{\ak{\thetaN}} (\boldsymbol x_m) \\
&\leq P_{\ak{\thetaN}} (\boldsymbol x_1) \\ 
&\leq 2^{N_{\mathcal{H}}} c(\ak{\thetaN}) \exp[ \max_{\boldsymbol h}f_{\ak{\thetaN}} (\boldsymbol x_1 , \boldsymbol h)]\\
&=2^{N_{\mathcal{H}}} c(\ak{\thetaN}) \exp[\min_{\boldsymbol x}\max_{\boldsymbol h}f_{\ak{\thetaN}} (\boldsymbol x , \boldsymbol h)]
\end{align*}
To prove Proposition \ref{prp:prop2}, we next expand the function \(f_{\ak{\thetaN}} (\boldsymbol x, \boldsymbol h)\) from \eqref{eq:f} as
\[
f_{\ak{\thetaN}} (\boldsymbol x, \boldsymbol h)= \sum_{j=1}^{N_\mathcal{H}} h_j \theta_{j,\ak{\thetaidx}}^{\mathcal{H}}  +   \sum_{i=1}^N \left( \theta_{i,\ak{\thetaidx}}^{\mathcal{V}} + \sum_{j=1}^{N_\mathcal{H}}  h_j  \theta_{ij,\ak{\thetaidx}}^{\mathcal{VH}}\right)x_i=\sum_{i=1}^N x_i  \theta_{i,\ak{\thetaidx}}^{\mathcal{V}} +  \sum_{j=1}^{N_\mathcal{H}} \left( \theta_{j,\ak{\thetaidx}}^{\mathcal{H}}  + \sum_{i=1}^N  x_i  \theta_{ij,\ak{\thetaidx}}^{\mathcal{VH}}\right)h_j.
\]
By this and the fact that \(x_i,h_j\in\{\pm 1\}\), we then have
\begin{align}
\nonumber \max_{\boldsymbol x}  f_{\ak{\thetaN}} (\boldsymbol x, \boldsymbol h) =  \sum_{j=1}^{N_\mathcal{H}} h_j \theta_{j,\ak{\thetaidx}}^{\mathcal{H}}  +   a_{\ak{\thetaN}, \mathcal{H}} (\boldsymbol h), & \min_{\boldsymbol x}  f_{\ak{\thetaN}} (\boldsymbol x, \boldsymbol h) =  \sum_{j=1}^{N_\mathcal{H}} h_j \theta_{j,\ak{\thetaidx}}^{\mathcal{H}}  - a_{\ak{\thetaN}, \mathcal{H}} (\boldsymbol h),\\
\nonumber \max_{\boldsymbol h}  f_{\ak{\thetaN}} (\boldsymbol x, \boldsymbol h) = \sum_{i=1}^{N}x_i \theta_{i,\ak{\thetaidx}}^{\mathcal{V}} + b_{\ak{\thetaN}, \mathcal{V}} (\boldsymbol x), & \min_{\boldsymbol h}  f_{\ak{\thetaN}} (\boldsymbol x, \boldsymbol h) = \sum_{i=1}^{N}x_i \theta_{i,\ak{\thetaidx}}^{\mathcal{V}}  -   b_{\ak{\thetaN}, \mathcal{V}} (\boldsymbol x), \\ 
\label{eq:max}
a_{\ak{\thetaN}, \mathcal{H}} (\boldsymbol h) \equiv \sum_{i=1}^N \left| \theta_{i,\ak{\thetaidx}}^{\mathcal{V}} + \sum_{j=1}^{N_\mathcal{H}}  h_j  \theta_{ij,\ak{\thetaidx}}^{\mathcal{VH}}\right|, &  b_{\ak{\thetaN}, \mathcal{V}} (\boldsymbol x) \equiv \sum_{j=1}^{N_\mathcal{H}} \left| \theta_{j,\ak{\thetaidx}}^{\mathcal{H}}  + \sum_{i=1}^N  x_i  \theta_{ij,\ak{\thetaidx}}^{\mathcal{VH}}\right|,
\end{align}
where \(\boldsymbol h^T \ak{\thetaN}^{\mathcal{H}}=\sum_{j=1}^{N_\mathcal{H}} h_j \theta_{j,\ak{\thetaidx}}^{\mathcal{H}}\), \(\boldsymbol x^T \ak{\thetaN}^{\mathcal{V}}= \sum_{i=1}^{N}x_i \theta_{i,\ak{\thetaidx}}^{\mathcal{V}}\) and \(\Gam\equiv \max_{\boldsymbol h} a_{\ak{\thetaN}, \mathcal{H}} (\boldsymbol h)\).
From this, it follows that
\begin{align*}
\REP_{(\boldsymbol X, \boldsymbol H)}(\ak{\thetaN}) &=  \max_{\boldsymbol h}\max_{\boldsymbol x}f_{\ak{\thetaN}} (\boldsymbol x , \boldsymbol h) -  \min_{\boldsymbol h}\min_{\boldsymbol x}f_{\ak{\thetaN}} (\boldsymbol x , \boldsymbol h)\\
&= \max_{\boldsymbol h_1 }  \max_{\boldsymbol h_2 }\left[ (\boldsymbol h_1 - \boldsymbol h_2)^T  \ak{\thetaN}^{\mathcal{H}}  +   a_{\ak{\thetaN}, \mathcal{H}} (\boldsymbol h_1)  + a_{\ak{\thetaN}, \mathcal{H}} (\boldsymbol h_2)\right],
\end{align*}
which leads to the upper bound \(\REP_{(\boldsymbol X, \boldsymbol H)}(\ak{\thetaN}) \leq 2 \Gam + 2 | \ak{\thetaN}^{\mathcal{H}} |_1\). Then, taking \(\boldsymbol h_1=\boldsymbol h_2\) (i.e., before maximization) gives \(\REP_{(\boldsymbol X, \boldsymbol H)}(\ak{\thetaN}) \geq 2\Gam\) and taking \(\boldsymbol h_1=-\boldsymbol h_2\), such that \(\boldsymbol h_1^T \ak{\thetaN}^{\mathcal{H}} = |\ak{\thetaN}^{\mathcal{H}}|_1\), gives \(\REP_{(\boldsymbol X, \boldsymbol H)}(\ak{\thetaN}) \geq 2|\ak{\thetaN}^{\mathcal{H}}|_1\); this yields the lower bound \(\REP_{(\boldsymbol X, \boldsymbol H)}(\ak{\thetaN})\geq 2\max\{\Gam, |\ak{\thetaN}^{\mathcal{H}} |_1\}\).

We next consider \(\elt\) and, by \eqref{eq:max}, write
\begin{align*}
\elt &=  \max_{\boldsymbol h}\max_{\boldsymbol x}f_{\ak{\thetaN}} (\boldsymbol x , \boldsymbol h) -  \max_{\boldsymbol h}\min_{\boldsymbol x}f_{\ak{\thetaN}} (\boldsymbol x , \boldsymbol h)\\
&= \max_{\boldsymbol h_1} \min_{\boldsymbol h_2 }\left[ (\boldsymbol h_1 - \boldsymbol h_2)^T \ak{\thetaN}^{\mathcal{H}} + a_{\ak{\thetaN}, \mathcal{H}} (\boldsymbol h_1)  + a_{\ak{\thetaN}, \mathcal{H}} (\boldsymbol h_2)\right].
\end{align*}
Taking \(\boldsymbol h_1=\boldsymbol h_2\) and maximizing over both \(\boldsymbol h_1,\boldsymbol h_2\) produces the upper bound \(\elt \leq 2\Gam\). Then, using \((\boldsymbol h_1 - \boldsymbol h_2)^T \ak{\thetaN}^{\mathcal{H}} + a_{\ak{\thetaN}, \mathcal{H}} (\boldsymbol h_2) \geq - 2|\ak{\thetaN}^{\mathcal{H}} |_1\) and maximizing over \(\boldsymbol h_1\) gives \(\elt \geq \Gam- 2|\ak{\thetaN}^{\mathcal{H}} |_1\), while setting \(\boldsymbol h_1=\boldsymbol h_2^*\) for \(\boldsymbol h_2^*\) such that \(-(\boldsymbol h_2^*) ^T \ak{\thetaN}^{\mathcal{H}} + a_{\ak{\thetaN}, \mathcal{H}} (\boldsymbol h_2^*) = \min_{\boldsymbol h_2} [-\boldsymbol h_2^T \ak{\thetaN}^{\mathcal{H}} + a_{\ak{\thetaN}, \mathcal{H}} (\boldsymbol h_2)]\) gives \(\elt \geq 2 a_{\ak{\thetaN}, \mathcal{H}} (\boldsymbol h_2^*) \geq \Gamc \equiv \min_{\boldsymbol h} a_{ \ak{\thetaN}, \mathcal{H}} (\boldsymbol h)\). Finally, note that for any \(\boldsymbol h\), the triangle inequality gives
\begin{align*}
\Gam  \equiv \max_{\boldsymbol h_1}  a_{\ak{\thetaN}, \mathcal{H}} (\boldsymbol h_1) &\geq
[a_{\ak{\thetaN}, \mathcal{H}} (\boldsymbol h)  + a_{\ak{\thetaN}, \mathcal{H}} (-\boldsymbol h)]/2 \\
&=
2^{-1}\sum_{i=1}^{N } \left(\left| \theta_{i,\ak{\thetaidx}}^{\mathcal{V}}   + \sum_{j=1}^{N_{\mathcal{H}}} h_j \theta_{ij,\ak{\thetaidx}}^{\mathcal{VH}} \right|
+ \left| \theta_{i,\ak{\thetaidx}}^{\mathcal{V}}   - \sum_{j=1}^{N_{\mathcal{H}}} h_j \theta_{ij,\ak{\thetaidx}}^{\mathcal{VH}} \right|\right)\\
& \geq \sum_{i=1}^{N }  \left| \theta_{i,\ak{\thetaidx}}^{\mathcal{V}} \right| \equiv | \ak{\thetaN}^{\mathcal{V}}|_1.
\end{align*}
\hfill \(\Box\)

\textbf{Proof of Theorem \ref{thm:something}.} Let \(L_{\ak{\thetaN}}(\boldsymbol X) = \log[ P_{ \ak{\thetaN}}(\boldsymbol X)/ \min_{\boldsymbol y \in \mathcal{X}^N} P_{\ak{\thetaN}}(\boldsymbol y) ]/\REP(\ak{\thetaN})\), where again \(\boldsymbol X=(X_1, \dots, X_N)\) and \(\REP(\ak{\thetaN})= \log[\max_{\boldsymbol y \in \mathcal{X}^N} P_{\ak{\thetaN}}(\boldsymbol y)/ \min_{\boldsymbol y \in \mathcal{X}^N} P_{ \ak{\thetaN}}(\boldsymbol y) ]\). As \(L_{\ak{\thetaN}}(\boldsymbol X)\in[0,1]\), convergence of \(L_{ \ak{\thetaN}}(\boldsymbol X)\) to 1 in probability under \(P_{\ak{\thetaN}}\) is equivalent to convergence to \(1\) in expectation under \(P_{\ak{\thetaN}}\) (i.e., convergence in expectation implies probabilistic convergence by Markov's inequality while probabilistic convergence implies convergence in expectation by uniform integrability/boundedness).

For \(\epsilon \in(0,1)\), let \(\mathcal{M}_{\epsilon, \ak{\thetaN}}\) denote a modal set as in \eqref{eq:mode}. By Theorem \ref{thm:degenFOES}, \(P_{\ak{\thetaN}}\left( \boldsymbol X\in \mathcal{M}_{\epsilon, \ak{\thetaN}}\right) \rightarrow 1\) holds as \(N \rightarrow \infty\) and, by definition of \eqref{eq:mode}, \(\boldsymbol X\in \mathcal{M}_{\epsilon, \ak{\thetaN}}\) follows if and only if
\(1-L_{\ak{\thetaN}}(\boldsymbol X)<\epsilon\). Hence, \(L_{\ak{\thetaN}}(\boldsymbol X) \stackrel{p,\mathrm{E}}{\longrightarrow} 1\) holds under \(\ak{\thetaN}\) in Theorem \ref{thm:something}. The convergence \(L_{-\ak{\thetaN}}(\boldsymbol X) \stackrel{p,\mathrm{E}}{\longrightarrow} 1\) under \(-\ak{\thetaN}\) likewise follows from Corollary \ref{cor:sign}. \hfill \(\Box\)

\BeginKnitrBlock{remark}
{}Consider the exponential graph model from \eqref{eq:mod2} involving counts of edges, 2-stars and triangles with fixed parameters \(\thetaN = (\theta_{1},\theta_2,\theta_3)^\prime \in \mathbb{R}^3\) for \(N\geq 1\). This model is unstable whenever \(|\theta_2| + |\theta_3| \neq 0\). To see this, consider an even number \(n>2\) of nodes and let \(\boldsymbol x_0\) denote the data outcome in \(\mathcal{X}^N \equiv \{0,1\}^N\) with all \(N = {n \choose 2}\) edges being zero, let \(\boldsymbol x_1\) denote the outcome with all edges being 1, and let \(\boldsymbol x_2\) denote the edge configuration from dividing the nodes into two equal groups, with no edges within a group and all edges between the groups (so that no triangles exist in \(\boldsymbol x_2\)). Then, the \(N\)-scaled log-ratio \eqref{eq:elpr} for the exponential graph model \eqref{eq:mod2} can, by definition, be bounded below by
\begin{align*}
\frac{1}{N}\REP_N(\ak{\thetaN}) &\geq \max_{i=1,2}\frac{1}{N}
\left| \log\left[ \frac{P_{\ak{\thetaN}}(\boldsymbol x_i)}{P_{\ak{\thetaN}}(\boldsymbol x_0)}\right] \right| \\
&= (n-2) \max\left\{ \left| \theta_2 + \frac{\theta_3}{3}+\frac{\theta_1}{n-2} \right|, \frac{n}{4(n-1)} \left| \theta_2 + \frac{8\theta_1}{n-2} \right| \right\};
\end{align*}
a similar expression also holds for an odd node number \(n>2\). Consequently, for all fixed parameters excluding \(\theta_2=\theta_3=0\), \(\lim_{N\to \infty}\REP_N(\ak{\thetaN})/N=\infty\) then follows and the graph model with 2-stars and triangles is S-unstable, as suggested by the breach of Proposition \ref{prp:prop1} for this model when \(|\theta_2|+|\theta_3|\neq 0\).
\EndKnitrBlock{remark}

\BeginKnitrBlock{remark}
{}Let \(M \geq 1\) denote a possible number of replications and consider data \(\boldsymbol Y_{N,M} \equiv (\boldsymbol X^{(1)}_N, \dots, \boldsymbol X^{(M)}_N)\) formed by \(\{ \boldsymbol X^{(j)}_N\}_{j=1}^M\) as \(M\) iid replications of a random vector \(\boldsymbol X_N=(X_1,\ldots,X_N)\), where the latter follows a FOES model with probabilities \(P_{\ak{\thetaN}}(\boldsymbol x)>0\), \(\boldsymbol x\in\mathcal{X}^N\). This leads to a joint model, say \(P_{\ak{\thetaN}}(\boldsymbol y)\), \(\boldsymbol y\in\mathcal{X}^{NM}\), for \(\boldsymbol Y_{N,M}\) consisting of \(N*M\) random variables in total. Then, the LREP for \(\boldsymbol Y_{N,M}\), scaled by associated size, is given by
\begin{align*}
\frac{1}{NM}\REP_{\boldsymbol Y_{N,M}}(\ak{\thetaN} ) &\equiv \frac{1}{NM}\log\left[\frac{\max_{\boldsymbol y \in \mathcal{X}^{NM}}P_{\ak{\thetaN}}(\boldsymbol y)  }{\min_{\boldsymbol y \in \mathcal{X}^{NM}}P_{\ak{\thetaN}}(\boldsymbol y)} \right] \\
&= \frac{1}{NM}\log\left[\frac{\max_{\boldsymbol x \in \mathcal{X}^{N}}P_{\ak{\thetaN}}(\boldsymbol x)  }{\min_{\boldsymbol x \in \mathcal{X}^{N}}P_{\ak{\thetaN}}(\boldsymbol x)} \right]^M \equiv \frac{1}{N}\REP_{\boldsymbol X_{N}}(\ak{\thetaN} ),
\end{align*}
where \(\REP_{\boldsymbol X_{N}}(\ak{\thetaN} ) \equiv \REP(\ak{\thetaN} )\) denotes the log-ratio of extremal probabilities for \(\boldsymbol X_N\) defined from \eqref{eq:elpr}. That is, due to iid properties, the sample-size corrected LREP for \(\boldsymbol Y_{N,M}\) equals the analog, \(\REP(\ak{\thetaN} )/N\), from the underlying common data model for \(\boldsymbol X_N\) alone, regardless of the level \(M \geq 1\) of independent replication. Hence, the definition of an S-unstable model is unaffected by independent replication. For computational purposes, this aspect also implies that if the original data \(\boldsymbol X_N=(X_1,\ldots,X_N)\) in a FOES model consist of \(N\) iid random variables, then the size-scaled log-ratio \eqref{eq:elpr} may be calculated as
\[
\frac{1}{N}\REP( \ak{\thetaN} ) \equiv \frac{1}{N}\REP_{\boldsymbol X_{N}}( \ak{\thetaN} ) = \log\left[ \frac{\max_{ x \in \mathcal{X}}P_{\ak{\thetaN}}(X_1=x)} {\min_{ x \in \mathcal{X}}P_{\ak{\thetaN}}(X_1=x)} \right]
\]
based on the extremal probabilities of just one random variable \(X_1\).
\EndKnitrBlock{remark}

\hypertarget{details-on-centered-graph-example}{%
\section{Details on Centered Graph Example}\label{details-on-centered-graph-example}}

\ak{To illustrate centered model parameterizations and n examination of stability in these, consider the two-star model (\ref{eq:mod2}) for the $N = {n \choose 2}$ edges in simple graph with $n$ nodes and binary edge-variables $(X_1,\ldots,X_N)\in\{0,1\}^N$. Here a common or standard parameterization in (\ref{eq:mod2}) leads to a conditional probability of ``$1$" for an edge $i$ as
$$
\mathrm{logit}[  P_{\thetaN}(X_{i} = 1 | X_j=x_j, j\neq i)=  \theta_1 + \theta_2 \sum_{j \in \mathcal{N}_i } x_j,$$ 
based on summing other edge observations $x_j$ in a neighborhood $\mathcal{N}_i = \{{s}_j: {s}_i \cap {s}_j \neq \emptyset\}$ to edge $i$ (i.e., edges $j$,  marked by pairs of graph vertices ${s}_{j} = \{ v_{j_1}, v_{j_2} \}$, that share a common vertex with edge $i$ marked by the vertex pair ${s}_i =\{v_{i_1}, v_{i_2}\}$).  In contrast, a centered conditional would yield 
$$
\mathrm{logit}[  P_{\thetaN}(X_{i} = 1 | X_j=x_j, j\neq i)=  \theta_1 +  \frac{\theta_2}{2(n-1)}\sum_{j \in \mathcal{N}_i} (x_{j}-\kappa),
$$
involving a parameter $\theta_1 \equiv \mathrm{logit}(\kappa)$ for $\kappa \in (0,1)$ and $2(n-2)$ as the size of $\mathcal{N}_i$ (cf. \citep{kaiser2009exploring}); the corresponding joint model would involve parameters $\theta_{1,\thetaidx}= \mathrm{logit}(\kappa) - \kappa\theta_2$
and $\theta_{2,\thetaidx} = \theta_2/(2(n-2))$ in (\ref{eq:mod2}). The purpose of the centerization is to have $\kappa \in(0,1)$ represent a model mean parameter (note $E X_i =\kappa$ is the edge proportion/probability under independence $\theta_2=0$), while a separate parameter $\theta_2$ (for dependence) modifies the conditional probability of ``$X_i=1$" up/down from $\kappa$, depending on neighbors $x_j=1$ or 0.  A similar interpretation does not hold in the uncentered model; see \citep{kaiser2009exploring} and \citep{casleton2017local} for a discussion of the centered parameterization in spatial and network modeling, where the intent is to improve parameter interpretation and help detect degeneracy (e.g., intuitively given by large dependence parameters $\theta_2$, so that $\kappa$ fails to correspond to the mean or the marginal proportion of $1$'s in data generations from the model).  However, from a different perspective using the model measures here, centered parameterizations also lead to models which are more stable over wider regions of the parameter space. To illustrate, in the standard/uncentered parameterization for a graph model (\ref{eq:mod2}) with only two-stars $\thetaN = (0,\theta_2,0)$, our measure $\REP(\thetaN)/N$ of model instability becomes $\REP(\thetaN)/N = |\theta_2| (n-2)$, which is unbounded as $N\to \infty$ for all non-zero parameters $\theta_2$ (i.e., all models with $\theta_2\neq 0$ are S-unstable).  However, in the centered parameterization  with only two-stars, corresponding to $\kappa=1/2$ and $\thetaN= 0.5*(-\theta_2,\theta_2/(n-2),0)$, the measure becomes
\begin{align*}
   \REP(\thetaN)/N  &= 0.5*|\theta_2|  \max_{x_i, y_i \in\{0,1\}}\left|  \frac{1}{N}\sum_{i=1}^N x_i \left( 1 - \frac{1}{2(n-2)} \sum_{j \in \mathcal{N}_i} x_j\right) - \frac{1}{N}\sum_{i=1}^N y_i \left( 1 - \frac{1}{2(n-2)} \sum_{j \in \mathcal{N}_i} y_j\right)\right| \\ 
   &\leq |\theta_2|,   
\end{align*}
which is bounded, for any fixed $\theta_2 \in\mathbb{R}$, as $N$ increases (i.e., note $\sum_{j \in \mathcal{N}_j} x_i/(2(n-2)) \in [0,1]$ is a conditional/neighborhood sample proportion while $\sum_{i=1}^N x_i/N \in[0,1]$ is a marginal proportion). This aspect owes to adjusting parameters by neighborhood sizes in centered conditional distributions, but centering also induces an additional effect of alternating signs in parameters (e.g., $\thetaN= 0.5*(-\theta_2,\theta_2/(n-2),0)$). The latter has been suggested in other contexts with exponential graph models for regulating degeneracy \citep[cf.][]{snijders2006new}.}

\clearpage

\bibliographystyle{imamat}
\bibliography{rbm.bib}

\end{document}